\documentclass[lettersize,journal]{IEEEtran}
\usepackage{cite}
\usepackage{amsmath,amsfonts}
\usepackage[linesnumbered,ruled,vlined]{algorithm2e}
\usepackage{array}
\usepackage[justification=centering]{caption}
\usepackage{textcomp}
\usepackage{stfloats}
\usepackage{url}
\usepackage{verbatim}
\usepackage{graphicx}
\usepackage{subfig}
\usepackage{tabularx}
\usepackage{hyperref}
\usepackage{booktabs}
\usepackage{diagbox}
\usepackage{makecell}
\usepackage{color}
\usepackage{fancyhdr}
\newtheorem{remark}{Remark}
\newtheorem{assumption}{Assumption}
\newtheorem{proposition}{Proposition}
\newtheorem{corollary}{Corollary}

\newtheorem{lemma}{Lemma}
\newtheorem{theorem}{Theorem}
\newtheorem{proof}{Proof}

\fancypagestyle{ieee_notice}{
    \fancyhf{}  
    \fancyfoot[C]{\footnotesize \textcolor{blue}{This work has been accepted for publication in IEEE Transactions on Systems, Man, and Cybernetics: Systems, available at: https://doi.org/10.1109/TSMC.2025.3565568. Copyright Statement: 2168-2216 © 2025 IEEE. All rights reserved, including rights for text and data mining, and training of artiﬁcial intelligence and similar technologies. Personal use is permitted, but republication/redistribution requires IEEE permission.}}
}

\begin{document}
\allowdisplaybreaks[4]
\title{\textit{Prox-DBRO-VR}:
A Unified Analysis on Byzantine-Resilient Decentralized Stochastic Composite Optimization with Variance Reduction}

\author{Jinhui Hu, Guo Chen, Huaqing Li, Xiaoyu Guo, Liang Ran, and Tingwen Huang
\thanks{The work is supported in part by the Fundamental Research Funds for the Central Universities of Central South University under grant 2023ZZTS0355, in part by the Fundamental Research Funds for the Central Universities under Grant SWU-XDJH202312, in part by the National Natural Science Foundation of China under Grant 62173278, and in part by the Chongqing Science Fund for Distinguished Young Scholars under Grant 2024NSCQ-JQX0103. \textit{(Corresponding author: Guo Chen)}.}
\thanks{Jinhui Hu is with the School of Automation, Central South University, Changsha 410083, Hunan, China, and also with the Department of Mechanical Engineering, City University of Hong Kong, Kowloon Tong, Kowloon, Hong Kong SAR, China (e-mail: jinhuihu3-c@my.cityu.edu.hk).}
\thanks{Guo Chen is with the School of Electrical Engineering and Telecommunications, University of New South Wales, Sydney, NSW 2052, Australia (e-mail: guo.chen@unsw.edu.au).}
\thanks{Huaqing Li and Liang Ran are with the Chongqing Key Laboratory of Nonlinear Circuits and Intelligent Information Processing, College of Electronic and Information Engineering, Southwest University, Chongqing 400715, China (e-mail: huaqingli@swu.edu.cn; ranliang\_rl@163.com).}
\thanks{Xiaoyu Guo is with the Department of Mechanical Engineering, City University of Hong Kong, Kowloon Tong, Kowloon, Hong Kong SAR, China (e-mail: xiaoyguo@cityu.edu.hk).}
\thanks{Tingwen Huang is with the Faculty of Computer Science and Control Engineering, Shenzhen University of Advanced Technology, Shenzhen 518055, China (e-mail: tingwen.huang@siat.ac.cn).}
}
\markboth{}%
{Shell \MakeLowercase{\textit{et al.}}: A Sample Article Using IEEEtran.cls for IEEE Journals}


\maketitle
\thispagestyle{ieee_notice}

\begin{abstract}
Decentralized stochastic gradient algorithms efficiently solve large-scale finite-sum optimization problems when all agents in the network are reliable. However, most of these algorithms are not resilient to adverse conditions, such as malfunctioning agents, software bugs, and cyber attacks. This paper aims to handle a class of general composite optimization problems over multi-agent systems (MASs) in the presence of an unknown number of Byzantine agents. Building on a resilient aggregation mechanism and the proximal-gradient mapping method, a Byzantine-resilient decentralized stochastic proximal-gradient algorithmic framework is proposed, dubbed \textit{Prox-DBRO-VR}, which achieves an optimization and control goal using only local computations and communications. To asymptotically reduce the noise variance arising from local gradient estimation and accelerate the convergence, we incorporate two localized variance-reduced (VR) techniques (\textit{SAGA} and \textit{LSVRG}) into \textit{Prox-DBRO-VR} to design \textit{Prox-DBRO-SAGA} and \textit{Prox-DBRO-LSVRG}. By analyzing the contraction relationships among the gradient-learning error, resilient consensus condition, and convergence error in a unified theoretical framework, it is proved that both \textit{Prox-DBRO-SAGA} and \textit{Prox-DBRO-LSVRG}, with a well-designed constant (resp., decaying) step-size, converge linearly (resp., sub-linearly) inside an error ball around the optimal solution to the original problem under standard assumptions. A trade-off between convergence accuracy and Byzantine resilience in both linear and sub-linear cases is also characterized. In numerical experiments, the effectiveness and practicability of the proposed algorithms are manifested via resolving a decentralized sparse machine-learning problem under various Byzantine attacks.
\end{abstract}


\section{Introduction}\label{Intro}
\subsection{Literature Review}\label{LiteRev}
Decentralized optimization has garnered significant attention and achieved substantial algorithmic breakthroughs in the fields of machine learning \cite{Xin2022,Kun2020,Xin2023}, smart grids \cite{Zhai2022}, cooperative control \cite{Li2021b}, and uncooperative games \cite{Huang2022}. Decentralized algorithms have advantages of high-efficiency in massive-scale optimization problems, good scalability over geographically distributed MASs, and a lower communication burden for the master/central agent in a parameter-server distributed structure. 

The increasing deployment of MASs in aforementioned mission-critical applications, there are unavoidable security issues in the course of optimization and control, such as poisoning data, software bugs, malfunctioning devices, cyber attacks, and privacy leakage. Many endeavors have been devoted to develop various privacy-preserving methods in decentralized optimization, for instance, differential privacy \cite{Wang2024, Xing2024}. Nevertheless, privacy issues are beyond the scope of this paper, as they neither directly cause compromised agents nor violate the prescribed update rules of algorithms. Other issues, however, may lead to node-level failures during multi-agent optimization and control, which is known as Byzantine problems \cite{Lamport1982}, where the malfunctioning or compromised agents are referred to as Byzantine agents. Byzantine agents colluding with each other are able to impede some notable Byzantine-free (non-resilient) decentralized optimization algorithms \cite{Saadatniaki2020,Xin2020f,Kun2020,Li2021b,Li2021,Pu2021,Qureshi2021,Huang2022,Li2022,Zhai2022,Xin2022,Wang2023} from achieving convergence, or even cause disagreement and divergence \cite{Sundaram2019,Yan2021}. For example, if a reliable agent is attacked and controlled by adversaries, the attacker can manipulate the agent to send misleadingly falsified information to its different reliable neighboring agents at each iteration. This can significantly disrupt the normal update of reliable agents, and even stop them from reaching consensus if a Byzantine agent deliberately sends misleading messages to their different reliable neighbors \cite{Yuan}. Therefore, researchers have been concentrating on designing resilient decentralized algorithms \cite{Yemini2022a,Zuo2022,El-Mhamdi2022,Fang2022,Li2022a,Wang2022,Wu2022,Peng2021,Li2019k} to alleviate or counteract the negative impact caused by Byzantine agents. In fact, there are various approaches to secure resilience in a decentralized manner. One popular line is to combine various screening or filtration techniques with decentralized optimization algorithms. To name a few, a pioneering work \cite{Sundaram2019} achieves Byzantine resilience via adopting trimmed mean (TM) aggregation, which discards a subset of the largest and smallest messages in the aggregation step. However, \cite{Sundaram2019} focuses on a class of scalar-valued optimization problems for decentralized optimization. A follow-up method \textit{ByRDiE} \cite{Yang2019c} combines decentralized coordinate gradient descent with TM aggregation to tackle the vector-valued optimization problem in decentralized learning. One imperfection of \textit{ByRDiE} stems from low efficiency in handling large-scale finite-sum optimization problems due to the implementation of one-coordinate-at-one-iteration update. Hence, \textit{BRIDGE} \cite{Fang2022} combines respectively four screening techniques including coordinate-wise trimmed-mean, coordinate-wise median, Krum function, and a combination of Krum and coordinate-wise trimmed mean, with decentralized gradient descent (\textit{DGD}) \cite{Nedic2009}. However, these four screening mechanisms either incur prohibitive computational overhead or impose extra restrictions on the number of neighbors or the network topology. To address these issues, building upon a centered clipping (\textit{CC}) technique \cite{He2022}, He et al. \cite{Praneeth2021} propose self-centered clipping (\textit{SCC}) for decentralized implementations,
which not only ensures Byzantine resilience but also provides a complete solution framework for decentralized non-convex optimization problems. One imperfection of \textit{SCC} \cite{He2022} is that the theoretical choice of localized clipping parameters requires knowledge of information that is assumed to be inaccessible. This limitation mandates manual parameter calibration in the practical deployment of \textit{SCC}, thereby compromising its autonomy. Another work \cite{Guo2021} designs a two-stage filtering technique to mitigate Byzantine threats, which is independent of any clairvoyant knowledge about Byzantine agents. Recent work \cite{Wu2022} systematically analyzes the relationship between two critical points, i.e., doubly-stochastic weight matrix and consensus, in the development of Byzantine-resilient decentralized optimization algorithms. On top of that, an iterative screening-based resilient aggregation rule, dubbed (\textit{IOS}), is designed in \cite{Wu2022}, which achieves Byzantine resilience and a controllable convergence error relied on the assumptions of bounded inner (node-level noise-based stochastic gradients) and outer (network-level aggregated gradients) variations.

Decentralized optimization algorithms \cite{Yang2019c,Guo2021,Fang2022,Wu2022,Wang2022} achieve Byzantine resilience via adopting various screening or filtering techniques. Nevertheless, the screening- or filtration-based methods may either impose a restriction on the number of neighboring agents or introduce additional computational complexity of $\mathcal{O}\left( mn \right)$ ($\mathcal{O}$, $m$, and $n$ denote an upper bound on the growth rate of the complexity, the total number of agents including both reliable and Byzantine agents in the network, and the dimension of decision variables, respectively), at each iteration (see \cite[TABLE II]{Fang2022}). This could be prohibitively expensive when the MAS is large-scale or the decision variable is high-dimensional. By integrating a resilient aggregation mechanism based on norm-penalized approximation \cite{Ben-ameur2016}, \textit{RSA} \cite{Li2019k} is developed for distributed federated-learning problems in the presence of Byzantine agents, which achieves Byzantine resilience without introducing prohibitive computational overhead. \cite{Peng2021} is a decentralized extension of \textit{RSA} \cite{Li2019k}. Via introducing a noise-shuffle strategy, \cite{Ma2022} enables \cite{Li2019k} to be differentially private
in distributed federated-learning tasks. However, the existing theoretical evaluation reveals that the algorithms proposed in \cite{Li2019k,Peng2021,Ma2022} sacrifice either convergence accuracy or convergence speed due to the employment of noise-based stochastic gradients. This dilemma can be resolved via employing VR techniques. Therefore, a recent Byzantine-resilient decentralized stochastic optimization algorithm \textit{DECEMBER} \cite{Peng2022} accelerates the convergence via incorporating VR techniques. Although \textit{DECEMBER} resolves the dilemma, it is still confined to empirical risk minimization (ERM) problems with local smooth objectives.

\subsection{Motivations}\label{Motiv}
On the one hand, Byzantine-resilient decentralized optimization algorithms \cite{Ben-ameur2016,Li2019k,Peng2021,Praneeth2021,Guo2021,Fang2022,Wang2022,Wu2022,Li2022a,Peng2022,He2022,Ma2022} are not available to handling optimization problems with a nonsmooth objective function, which is indispensable in many practical applications, such as sparse machine learning \cite{Ye2020,Alghunaim2021}, model predictive control \cite{Li2021b}, and signal processing \cite{Notarnicola2017a}. On the other hand, despite the fact that there are various decentralized optimization algorithms \cite{Notarnicola2017a,Li2019m,Ye2020,Alghunaim2021,Li2021,Xu2021,Xin2023} providing many insights to resolve the composite optimization problem, they all fail to consider any possible security issues over MASs. This renders the reliable agents under the algorithmic framework of \cite{Notarnicola2017a,Li2019m,Ye2020,Xu2021,Alghunaim2021,Li2021,Xin2023} vulnerable to Byzantine attacks or failures. To bridge this gap, this paper studies a category of composite optimization problems in the presence of Byzantine agents, where the local objective function associated with each agent consists of both smooth and nonsmooth objectives. In a nutshell, the non-trivial study on Byzantine-resilient decentralized stochastic composite optimization deserves further investigation, which features the main motivation of this paper. The integration of VR techniques serves as a side motivation, which not only helps reduce the prohibitive per-iteration gradient-computation cost when evaluating the local batch gradients \cite{Notarnicola2017a,Yang2019c,Sundaram2019,Li2019m,Alghunaim2021,Xu2021,Fang2022,Wang2022,Yemini2022a} but also asymptotically reduces the noise variance arising from local gradient estimation \cite{Li2019k,Praneeth2021,Guo2021,Peng2021,Qureshi2021,He2022,Wu2022,Li2022a}. In contrast to a recent notable algorithm (event-triggered) \textit{MW-MSR} \cite{Yuan,Yuan2023a}, which proposes an alternative approach to achieving Byzantine-resilient consensus with asynchronous multi-hop communication through detection and filtration, this paper focuses on decentralized optimization problems. In such problems, the heterogeneity in local gradient evaluations can lead to greater variance among reliable agents' states \cite{He2022}, which may cause false filtration when the reliable information is mistakenly identified as Byzantine information. This observation motivates the design of a screening- and filtration-free Byzantine-resilient decentralized optimization algorithm tailored to handle the worst-case scenario of Byzantine problems in decentralized optimization. The worst-case scenario indicates that the number of Byzantine agents is unknown and they can be omniscient to send any misleadingly falsified messages to their reliable neighbors.

\subsection{Contributions}\label{Contri}
\begin{enumerate}
\item This paper develops a Byzantine-resilient decentralized stochastic proximal-gradient algorithmic framework, dubbed \textit{Prox-DBRO-VR}, to resolve a class of composite (smooth + nonsmooth) finite-sum optimization problems over MASs under the worst case. The challenge of studying the nonsmooth objective function in the presence of Byzantine agents stems from incurring additional error terms consisting of coupled Byzantine and reliable information in contrast to Byzantine-free composite optimization \cite{Notarnicola2017a,Li2019m,Ye2020,Alghunaim2021,Li2021,Xu2021,Xin2023} and Byzantine-resilient smooth optimization \cite{Ben-ameur2016,Li2019k,Peng2021,Praneeth2021,Guo2021,Fang2022,Wang2022,Wu2022,Li2022a,Peng2022,He2022,Ma2022}. To handle these introduced error terms, we explore and seek the upper bounds on the subdifferentials related to Byzantine agents after applying the non-expansiveness of the proximal operator in the linear convergence case and bounded-gradient condition in the sub-linear convergence case, respectively.

\item Inspired by \cite{Ye2020,Xin2020f}, we incorporate localized versions of two VR techniques \textit{SAGA} \cite{Defazio2014c} and \textit{LSVRG} \cite{Kovalev2019}, into \textit{Prox-DBRO-VR}, to propose two Byzantine-resilient decentralized stochastic proximal gradient algorithms, namely \textit{Prox-DBRO-SAGA} and \textit{Prox-DBRO-LSVRG}, both of which reduce the per-iteration gradient-computation cost in evaluating batch gradients \cite{Notarnicola2017a,Yang2019c,Sundaram2019,Li2019m,Alghunaim2021,Xu2021,Fang2022,Wang2022,Yemini2022a} and circumvent the bounded-variance assumption required by noise-based stochastic gradient methods \cite{Li2019k,Praneeth2021,Guo2021,Peng2021,Qureshi2021,He2022,Wu2022,Li2022a,Ma2022}. The challenge of studying VR techniques in the presence of Byzantine agents lies in seeking an appropriate (Lyapunov) candidate function with respect to (w.r.t.) the gradient-learning sequences and errors, we address this challenge by exploiting the Bregman divergence, choosing a proper decaying or constant step-size, and utilizing appropriate intermediate constants (see the proofs of Theorems \ref{T2}-\ref{T3}).

\item  The proposed algorithmic framework \textit{Prox-DBRO-VR} achieves Byzantine resilience without incurring prohibitive computational overhead in contrast to Byzantine-resilient decentralized optimization algorithms \cite{Sundaram2019,Yang2019c,Fang2022,Wang2022,Guo2021,Wu2022,Li2022a,Yuan} introducing additional computational complexity of $\mathcal{O}\left( mn \right)$ for screening or filtration processes. The theoretical analysis of \textit{Prox-DBRO-SAGA} and \textit{Prox-DBRO-LSVRG} imposes no assumption or restrictions on the number or proportion of Byzantine agents in the network, while only assumes a connected network among all reliable agents (see Corollary \ref{C1}), which is less restrictive than the topology condition of many related studies, such as \cite{Sundaram2019,Yang2019c,Fang2022,Wang2022,Wu2022}. Under this assumption, theoretical results reveal an explicit trade-off between convergence accuracy and Byzantine resilience in both cases (see Theorems \ref{T2}-\ref{T3} for details), providing directions to optimize the performance of the proposed algorithms in practice.
\end{enumerate}

\subsection{Organization}\label{Organ}
We provide the remainder of the paper in this part. Section \ref{Preli} presents the basic notation, problem statement, problem reformulation, and setup of its robust variant. The connection of the proposed algorithms with existing methods and the algorithm development are elaborated in Section \ref{AlgoDeve}. Section \ref{ConAna} details the convergence results of the proposed algorithms. Case studies on decentralized learning problems with various Byzantine attacks to illustrate the effectiveness and performance of the proposed algorithms are carried out in Section \ref{ExpRe}. Section \ref{Conclu} concludes the paper and states our future direction. Some detailed derivations are placed to Appendix for coherence.

\section{Preliminaries}\label{Preli}
\subsection{Basic Notation}\label{BasNo}
Throughout the paper, we assume all vectors are column vectors if there is no other specified. For arbitrary three vectors $\tilde x, \tilde y, \tilde z\in \mathbb{R}^{ n}$,
a positive scalar $a$ and a closed, proper, convex function $g:{\mathbb{R}^n} \to \mathbb{R} $, the proximal operator of $g$ is defined by
${\mathbf{prox}}_{{a, g}}\left \{ \tilde x \right \} = \arg {\min _{\tilde y \in \mathbb{R}^{n}}} \left \{ {g\left( \tilde y \right) + \frac{1}{{2a}}{{\left\| {\tilde y - \tilde x} \right\|}_2^2}} \right \}$;
let $\partial g\left( \tilde x \right)$ denote the subdifferential of $g$ at $\tilde x$, such that
\begin{equation*}
\partial g\left( \tilde x \right) = \left\{ { \tilde y |  \forall \tilde z \in {\mathbb{R}^n}, g\left(\tilde x \right) + \left\langle {\tilde y, \tilde z - \tilde x} \right\rangle  \le g\left( \tilde z \right)} \right\},
\end{equation*}
let ${{\partial _{\tilde x}}g\left( \tilde x \right)}$ denote one subgradient of $g$ at $\tilde x$, such that ${{\partial _{\tilde x}}g\left( \tilde x \right)} \in {{\partial}g\left( \tilde x \right)}$.
\begin{table}[!h]
\centering
\caption{Basic notations.}
\begin{tabularx}{8.8cm}{lX}  
\hline                      
\bf{Symbols}  & \bf{Definitions}  \\
\hline
${\mathbb R}$, ${{\mathbb R}^n}$, ${{\mathbb R}^{m \times n}}$ & the set of real numbers, $n$-dimensional column real vectors, $m \times n$ real matrices, respectively\\
${\mathbb{E}_{ \xi} }\left[ \cdot \right]$ & the expectation with respect to a random variable $ \xi $ \\
:= & the definition symbol\\
${I_n}$   & the $n \times n$ identity matrix\\
${0_{n}}$   & an $n$-dimensional column vector with all-zero elements\\
$1_m$   & an $m$-dimensional column vector with all-one elements\\
${\cdot ^ \top }$ & transpose of any matrices or vectors \\
$X \le Y$ & each element in $Y - X$  is nonnegative, where $X$ and $Y$ are two vectors or matrices with same dimensions\\
$\tilde x \otimes \tilde y$ & the Kronecker product of vectors $\tilde x$ and $\tilde y$\\
$\left| \cdot \right|$ & the operator to represent the absolute value of a constant or the cardinality of a set\\
$\left\|  \nu   \right\|_a$ & either the $a$-norm of $\nu \in {\mathbb{R}^n}$ equivalent to ${\left( {\sum\nolimits_{i = 1}^n {{{\left| {{\nu_i}} \right|}^a}} } \right)^{\frac{1}{a}}}$, $a \ge 1$, or its induced matrix norm.\\
${\mathbf{col}}\left\{ \cdot \right\}_i$ & the operator that vertically concatenates the input vectors indexed by $i$ into one column vector\\
${\lambda _{\min }}\left( X \right)$ & the minimum nonzero singular value of any matrix $X$\\
${\lambda _{\max }}\left( X \right)$ & the maximum singular value of any matrix $X$\\
\hline
\end{tabularx}
\label{Table 1}
\end{table}
The remaining basic notations of this paper are summarized in Table \ref{Table 1}.
\subsection{Problem Statement}\label{ProSta}
A network of $m$ agents connect with each other over an undirected network $\mathcal{G} := \left( {\mathcal{R} \cup \mathcal{B},\mathcal{E}} \right)$, where $\mathcal{R}$  ($2 \le \left| \mathcal{R} \right| \le m$) and $\mathcal{B}$ indicate the sets of reliable and Byzantine agents, respectively, and $\mathcal{E}$ represents the set of undirected (communication) edges among all agents. We represent the information flow from agent $j$ to its neighboring agent $i$ by the edge $(i, j) \in \mathcal{E}$, with the convention that $i < j$ always holds for $ i \in \mathcal{V}, j \in {\mathcal{N}_i}$. The mutual target of all reliable agents is to minimize (min) a general decentralized composite ERM problem as follows:
\begin{equation}\label{E2-1}
\mathop {\min }\limits_{\tilde x } \sum\limits_{i \in \mathcal{R}} {{f_i}\left( {\tilde x} \right) + {g}\left( {\tilde x} \right)},
\end{equation}
where $\tilde x \in {\mathbb{R}^n}$ is the decision variable, and ${f_i}:{\mathbb{R}^n} \to \mathbb{R}$ and ${g}:{\mathbb{R}^n} \to \mathbb{R}$, $i \in \mathcal{R}$, are two different objective functions. The local objective function $f_i$ can be further decomposed as ${f_i}\left( {\tilde x} \right) = \sum\nolimits_{l = 1}^{{q_i}} {f_i^l\left( {\tilde x} \right)}/q_i$, while the function $g$ serves as a shared nonsmooth objective among all reliable agents similar to literature \cite{Alghunaim2021,Ye2020,Xu2021}. Based on the above definitions, problem (\ref{E2-1}) is also recognized as a composite ERM problem in the presence of Byzantine agents. We assume that the optimal solution to (\ref{E2-1}) exists, denoted by $\tilde x^*$, and the local sample set associated with agent $i$ as ${\mathcal{Q}_i} = \left\{ {1, 2, \ldots ,{q_i}} \right\}$, $\forall i \in \mathcal{R}$. This paper aims to resolve the general composite ERM problem (\ref{E2-1}) under Byzantine attacks or failures, including but not limited to the zero-sum attack \cite{He2022}, Gaussian attack \cite{Wu2022}, and same-value attack \cite{Peng2021}, which will be testified in numerical experiments. In fact, failures, such as agent breakdown and possible disconnection of communication links, can be deemed as a category of not malicious Byzantine problems. We next specify the studied problem via the following standard assumptions.
\begin{assumption}(Convexity and smoothness).\label{Assu1}\\
a) For $i \in \mathcal{R}$, the local objective function $f_i$ is assumed to be $\mu _i$-strongly convex, and the local component objective function $f_i^l$ is $L_i^l$-smooth, $\forall l \in {\mathcal{Q}_i}$, i.e., $\forall {\tilde x}, {\tilde z} \in {\mathbb{R}^n}$,
\begin{subequations}\label{E2-2}
\begin{align}
&\label{E2-2-1}{\mu _i}\left\| {\tilde x - \tilde z} \right\|_2^2 \le {\left( {\nabla {f_i}\left( {\tilde x} \right) - \nabla {f_i}\left( {\tilde z} \right)} \right)^ \top }\left( {\tilde x - \tilde z} \right),\\
&\label{E2-2-2}{\left\| {\nabla f_i^l\left( {\tilde x} \right) - \nabla f_i^l\left( {\tilde z} \right)} \right\|_2} \le L_i^l{\left\| {\tilde x - \tilde z} \right\|_2},
\end{align}
\end{subequations}
where ${\mu _i}$ and $L_i^l$ are the strongly-convex and smooth parameters, respectively. \newline
b) The objective function $g$ is convex and not necessarily smooth.
\end{assumption}

For convenience, we define $\mu := \mathop {\min }\nolimits_{i \in \mathcal{R}} \left\{ {{\mu _i}} \right\}$ and $L: = {\max _{i \in \mathcal{R}}}\left\{ {{L_i}} \right\}$ with ${L_i}: = \sum\nolimits_{l = 1}^{{q_i}} {L_i^l} /{q_i}$. It can be verified that the global objective function is
$L$-smooth and $\mu$-strongly convex, with a condition number denoted by ${\kappa _f}: = L/\mu$.
\begin{remark}\label{R2-1}
Assumption \ref{Assu1}-a) is standard in recent literature \cite{Li2022,Pu2021,Li2021,Qureshi2021,Saadatniaki2020}. According to \cite[Chapter 3]{Bubeck2015c}, we know that $ 0< \mu \le L$, which indicates $\kappa_f \ge 1$. Moreover, in view of (\ref{E2-2-2}), it can be verified that the local objective functions $f_i$, $i \in \mathcal{R}$, are $L$-smooth as well. Under Assumption \ref{Assu1}, the optimal solution $\tilde x^*$ to (\ref{E2-1}) exists uniquely. The consideration of the possibly nonsmooth term $g$ is meaningful, which finds substantial applications in various fields, such as the standard 1-norm regularization term in sparse machine learning \cite{Notarnicola2017a,Ye2020,Alghunaim2021}, a nonsmooth indicator function in model predictive control \cite{Li2021b} to handle equality and set constraints, and a nonsmooth indicator function in energy resource coordination \cite{Li2021} to handle inequality and set constraints.
\end{remark}
\begin{assumption}(Network connectivity)\label{Assu2}
All reliable agents form a static network, denoted as $ {{\cal G}_{\cal R}} := \left( {\mathcal{R},{\mathcal{E}_\mathcal{R}}} \right)$, which is bidirectionally connected.
\end{assumption}
\begin{remark}\label{R2-2}
Assumption \ref{Assu2} is standard in recent literature \cite{Ye2024,Ma2022,He2022,Yemini2022a}, which implies that each reliable agent must have at least one reliable neighbor and an arbitrary number of Byzantine neighbors to enable the communication with any other reliable agents in the network. There are many examples satisfying Assumption \ref{R2-2} and a straightforward instance is the full-connected network among all agents. From another perspective, there must be two reliable agents that fail to exchange messages if ${{\cal G}_{\cal R}}$ is disconnected. In this scenario, the disconnected agent is also judged as a Byzantine agent subject to possible failures. Assumption \ref{Assu2} imposes few restrictions on the reliable network ${{\cal G}_{\cal R}}$ than the pioneering literature \cite{Sundaram2019}, which requires at least $\left| \mathcal{B} \right| +1$ paths between any two reliable agents when there are $\left| \mathcal{B} \right|$ Byzantine agents in the network. In view of this, there are some network examples, for instance the Dumbbell network \cite{Praneeth2021} satisfying Assumption \ref{Assu2} but violating the resilient network assumption made in \cite{Sundaram2019}.
\end{remark}
\subsection{Problem Reformulation}\label{ProRef}
To guarantee all reliable agents reach a consensus at the optimal solution, we need to reformulate (\ref{E2-1}) into an equivalent consensus problem. To achieve this goal, a global decision vector $x = {\mathbf{col}}{\left\{ {{x_i}} \right\}_{i \in \mathcal{R}}} \in {\mathbb{R}^{\left| \mathcal{R} \right|n}}$ containing $\left| \mathcal{R} \right|$ local copies of the decision variable $\tilde x$, is introduced, subject to (s.t.) the consensus constraint ${x_i} = {x_j},\left( {i,j} \right) \in \mathcal{E}_\mathcal{R}$. Therefore, it is natural to rewrite (\ref{E2-2}) as
\begin{equation}\label{E2-3}
\begin{aligned}
  &\mathop {\min }\limits_{x } F\left( x \right) + G\left( x \right), \hfill \\
  &{\text{s}}{\text{.t}}{\text{. }}{x_i} = {x_j},\left( {i,j} \right) \in {\mathcal{E}_\mathcal{R}}, \hfill \\
\end{aligned}
\end{equation}
where $F\left( x \right) := \sum\nolimits_{i \in \mathcal{R}} {f_i\left( {{x_i}} \right)}$ and $G\left( x \right) := \sum\nolimits_{i \in \mathcal{R}} {g\left( {{x_i}} \right)}$.

\subsection{Resilient Consensus Problem Setup}\label{ProRef}
To facilitate the design of the resilient aggregation mechanism, we examine a norm-penalized approximation variant of the consensus problem (\ref{E2-3}), originally studied in \cite{Ben-ameur2016}, as follows:
\begin{equation}\label{E2-4}
{x^*}: = \arg \mathop {\min}\limits_x \sum\limits_{i \in \mathcal{R}} {( {{f_i}\left( {{x_i}} \right) + g\left( {{x_i}} \right) + \frac{\phi }{2}\sum\limits_{j \in {\mathcal{R}_i}} {{{\left\| {{x_i} - {x_j}} \right\|}_a}} } )},
\end{equation}
where $a \ge 1$, $\phi$ is the penalty parameter, and ${\mathcal{R}_i}$ denotes the set of reliable neighbors of agent $i$. The norm penalty takes a soft-approximation replacement of the consensus constraint, thereby admitting a resilient distance between any two neighboring states. The distance is controlled by the penalty parameter $\phi$, indicating that a larger $\phi$ can bring a smaller gap between $x_i$ and $x_j$, $\left( {i,j} \right) \in {\mathcal{E}_\mathcal{R}}$. Therefore, the problem formulation (\ref{E2-4}) is regarded as a soft approximation of (\ref{E2-3}) as the former tolerates the state dissimilarities among any neighboring agents, which exhibits resilience to a prevalent issue of data silos in decentralized learning tasks. The equivalence between the soft-approximation problem (\ref{E2-4}) and the consensus problem (\ref{E2-3}) with respect to the original problem (\ref{E2-1}) is proved in Theorem \ref{T1}.

\section{Algorithm Development}\label{AlgoDeve}
\subsection{Connection to Existing Work}
Lian et al. in \cite{Lian2017} design a decentralized stochastic gradient descent algorithm, namely \textit{D-PSGD}, to resolve efficiently the consensus problem (\ref{E2-3}) in an ideal situation. The ideal situation fails to consider the presence of any malfunctioning or malicious agents, i.e., Byzantine agents, which may not be avoided in practical applications \cite{Li2022a,Fang2022,Ramanan2022,Wu2022,El-Mhamdi2022,Lamport1982,Zuo2022}. We next find out the reason why \textit{D-PSGD} cannot be applied directly to solving (\ref{E2-3}) when there are Byzantine agents in the network, and then seek out a feasible improvement to secure Byzantine resilience. We first recap the updates of the generalized \textit{D-PSGD} as follows:
\begin{subequations}\label{E3-0}
\begin{align}
&\label{E3-0-1}{{\bar x}_{i,k}} = {x_{i,k}} - {\alpha _k} {\nabla f_i\left( {{x_{i,k}}} \right)} , \\
&\label{E3-0-2}{x_{i,k + 1}} = \sum\limits_{j \in {{{\cal R}_i} \cup {{\mathcal{B}}_i}}} {{w_{ij}}{{{v_{ij,k}}}}},
\end{align}
\end{subequations}
where ${\mathcal{B}_i}$ denotes the set of Byzantine neighbors of agent $i$, $ i \in \mathcal{R}$.
${v_{ij,k}}: = \left\{ \begin{array}{l}
\!\!\!{{\bar x}_{j,k}},j \in {{\mathcal R}_i}\\
\!\!\!{z_{ij,k}},j \in {{\mathcal{B}}_i}
\end{array} \right.$ with ${z_{ij,k}}$ denoted as untrue or misleading information sent by Byzantine agent $j$, $j \in {\mathcal{B}}_i$, $\alpha_k$ denotes a proper constant or decaying step-size, ${\nabla f\left( {{x_{i,k}}} \right)}$ is the local batch gradient, ${w_{ij}}$ is the $i$-th row and $j$-th column element of a doubly stochastic weight matrix s.t. $\sum\nolimits_{j \in {{\mathcal N}_i}} {{w_{ij}}}  = \sum\nolimits_{j \in {{\mathcal N}_i}} {{w_{ji}}}  = 1$. Note that both $\mathcal{R}_i$ and $\mathcal{B}_i$ exclude agent $i$ (itself). If there is a reliable agent $i$ with a Byzantine neighbor $b$, $b \in {{\mathcal{B}}_i}$, then ${{z}_{ib,k}}$ could be untrue or misleading information (depending on whether agent $b$ experiences malfunction or manipulated by adversaries), to its reliable neighboring agents at $k$-th iteration. Since any Byzantine agent is assumed to be omniscient and able to learn from update rules, if agent $b$ sends an elaborately falsified message to its reliable neighbor $i$, then ${x_{i,k + 1}}$ could arbitrarily deviate from its true model.
For instance, Byzantine agent $b$ can blow $x_{i,k+1}$ up to infinity through continually transmitting a vector with infinite elements to its reliable neighbor $i$. Another example is that Byzantine agent $b$ can deter all reliable neighbors from achieving consensus at iteration $k$, via sending $\tilde x_{ib,k}$ with various values to its different reliable neighbors $i \in {\mathcal{R}_b}$. The primary cause of the aforementioned issues lies in the vulnerability of the aggregation step (\ref{E3-0-2}) to Byzantine agents. In fact, similar security threats also prevail in decentralized work \cite{Li2021,Saadatniaki2020,Xin2022,Pu2021,Qureshi2021,Li2022,Xin2020f,Kun2020}.
To address this issue, the \textit{SGD} family receives two important extensions, \textit{RSA} \cite{Li2019k} and \cite{Peng2021}, both of which achieve Byzantine resilience based on a resilient aggregation mechanism \cite{Ben-ameur2016}, where \cite{Peng2021} is a decentralized extension of \textit{RSA} \cite{Li2019k}. The theoretical analysis of both \textit{RSA} \cite{Li2019k} and \cite{Peng2021} is based on a bounded-variance assumption on the local stochastic gradient. With this assumption and the other standard assumptions (see \cite{Peng2021} for details), the sequence ${\left\{ {{x_k}} \right\}_{k \ge 0}}$ generated by the decentralized algorithm proposed in \cite{Peng2021} takes a convergent form of
\begin{equation}\label{E3-1}
\begin{aligned}
\mathbb{E}\left[ {\left\| {{x_{k + 1}} - {1_{\left| \mathcal{R} \right|}}\! \otimes \! {\tilde x^*}} \right\|_2^2} \right] \le  & \left( {1 - \eta {\alpha _k}} \right)\mathbb{E}\left[ {\left\| {{x_k} - {1_{\left| \mathcal{R} \right|}} \! \otimes \! {\tilde x^*}} \right\|_2^2} \right] \\
& + \alpha _k^2{\Delta _0} + {\alpha _k}{\Delta _1},
\end{aligned}
\end{equation}
where $\eta $ is a positive constant satisfying $0 < \eta {\alpha _k} < 1$, ${\Delta _0} := \sum\nolimits_{i \in \mathcal{R}} {32n{\phi ^2}{{\left| {{\mathcal{R}_i}} \right|}^2}}  + 4n{\phi ^2}{\left| {{\mathcal{B}_i}} \right|^2} + 2\sigma _i^2$ (${\sigma _i}>0$ is the bounded variance yielded by the biased estimation of the local batch gradients) and ${\Delta _1}: = \left( {n{\phi ^2}/\gamma } \right)\sum\nolimits_{i \in {\cal R}} {{{\left| {{{\cal B}_i}} \right|}^2}} $. Based on (\ref{E3-1}), one can establish either a sub-linear convergence rate with a smaller convergence error determined by the number of Byzantine agents, or a faster linear convergence rate with a larger convergence error determined jointly by the number of Byzantine agents and the bounded variance. In fact, this bounded variance ($\sigma _i^2$) exists commonly in recent literature, such as \cite{Praneeth2021,Guo2021,Peng2021,Qureshi2021,He2022,Wu2022,Li2022a}. Therefore, this paper aims to asymptotically reduce this bounded variance in the linear convergence case and circumvent the bounded-variance assumption. Inspired by the recent exploration of decentralized VR stochastic gradient algorithms \textit{diffusion-AVRG} \cite{Kun2020}, \textit{S-DIGing} \cite{Li2022}, \textit{GT-SAGA}/\textit{GT-SVRG} \cite{Xin2020f}, and \textit{GT-SARAH} \cite{Xin2022} that seek the solution to an optimization problem under an ideal Byzantine-free situation, we introduce two popular localized variance-reduced techniques \textit{SAGA} \cite{Defazio2014c} and \textit{LSVRG} \cite{Kovalev2019} to asymptotically reduce the noise variance in local gradient estimation. These two VR techniques allow us to derive a unified theoretical result on Byzantine-resilient decentralized stochastic proximal-gradient optimization.

\subsection{A General Algorithmic Framework}
Based on the above analysis, we propose a Byzantine-resilient decentralized stochastic-gradient algorithmic framework in Algorithm \ref{Algo1} to resolve (\ref{E2-4}) in the presence of Byzantine agents. Note that we denote temporarily the local stochastic gradient by $r_{i,k}$, which will be specified in Steps 4-5 of Algorithms \ref{Algo2}-\ref{Algo3}.

\begin{algorithm}[!h]
	\small
    \SetKwInput{KwInit}{Initialize}
	\SetKwBlock{Repeat}{Repeat for $ k = 0,1,2, \dots$}{}
	\SetKwBlock{End}{End for a required criterion is satisfied}{}
	\KwIn{a proper constant or decaying step-size $\alpha_k>0$.}
	\KwInit{an arbitrary starting point $x_{i,0} \in {\mathbb{R}^n}$ and a proper penalty parameter ${\phi}>0$.}
	\For{$k=0,1,\ldots,$}{
		\For{\rm{each reliable agent} $i \in \mathcal{R}$}{
			{\textbf{Transmit} its current local model ${x_{i,k}}$ to its neighbors $j \in {\mathcal{N}_i}$ and receive the true information ${x_{j,k}}$ or untrue information ${z_{ij,k}}$ from its neighbors;}\\
			{\textbf{Evaluate} the local stochastic gradient $r_{i,k}$;}\\
			{\textbf{Calculate} an intermediate variable according to the local resilient stochastic subgradient descent step:}
            \begin{equation*}\label{E3-2}
            {\bar x_{i,k}} = {x_{i,k}} - {\alpha _k}r_{i,k}- {\alpha _k}{\phi}\sum\limits_{j \in {\mathcal{N}_i}} {{\partial _{{x_i}}}{{\left\| {{x_{i,k}} - v_{ij,k}} \right\|}_a}},
            \end{equation*}
            where ${v_{ij,k}} := \left\{ \begin{gathered}
              {x_{j,k}}, {\text{if}} \; j \in {\mathcal{R}_i} \hfill \\
              {z_{ij,k}}, {\text{if}} \; j \in {\mathcal{B}_i} \hfill \\
            \end{gathered}  \right..$\\
			{\textbf{Update} its current local model according to the local proximal gradient step:}
                \begin{equation*}\label{E3-2+}
                {x_{i,k + 1}} = \arg \mathop {\min }\limits_{\tilde x \in {\mathbb{R}^n}} \left\{ {g\left( {\tilde x} \right) + \frac{1}{{2{\alpha _k}}}\left\| {\tilde x - {{\bar x}_{i,k}}} \right\|_2^2} \right\}.
                \end{equation*}

		}
		\For{\rm{each Byzantine agent} $i \in \mathcal{B}$}{\textbf{Send} an arbitrary vector $ {z_{ij,k}}$ to its neighbor $j$, $j \in {\mathcal{N}_i}$.}
	}
	\KwOut{all decision variables ${x_{i,k}}$, $i \in \mathcal{V}$, until a prescribed criterion is satisfied.}
	\caption{\textit{Prox-DBRO-VR} Framework.}
    \label{Algo1}
\end{algorithm}

\begin{remark}\label{R1}
Inspired by \cite{Ben-ameur2016}, the Byzantine resilience of \textit{Prox-DBRO-VR} is achieved by employing the resilient aggregation based on norm-penalized approximation. The literature \cite{Li2019k,Ma2022} leverages this strategy to handling distributed federated learning problems, and \cite{Peng2021} studies it in a decentralized manner. However, all these algorithms \cite{Li2019k,Peng2021,Ma2022} not only rely on a bounded-variance assumption in theoretical analysis but also exhibit a fundamental trade-off between convergence accuracy and speed. Therefore, one important motivation behind the design of \textit{Prox-DBRO-VR} is to bypass the bounded-variance assumption while resolving the fundamental trade-off, which can be achieved with the aid of VR techniques.
\end{remark}
\subsection{\textit{Prox-DBRO-SAGA} and \textit{Prox-DBRO-LSVRG}}
We introduce the localized version of two popular centralized VR techniques \textit{SAGA} \cite{Defazio2014c} and \textit{LSVRG} \cite{Kovalev2019}, into \textit{Prox-DBRO-VR}, to develop \textit{Prox-DBRO-SAGA} and \textit{Prox-DBRO-LSVRG}. The detailed updates of \textit{Prox-DBRO-SAGA} and \textit{Prox-DBRO-LSVRG} are presented in Algorithms \ref{Algo2}-\ref{Algo3}, respectively.
\begin{algorithm}[!htp]
	\small
    \SetKwInput{KwInit}{Initialize}
	\SetKwBlock{Repeat}{Repeat for $ k = 0,1,2, \dots$}{}
	\SetKwBlock{End}{End for a required criterion is satisfied}{}
	\KwIn{a proper constant or decaying step-size $\alpha_k>0$.}
	\KwInit{the same parameters and starting points according to Algorithm \ref{Algo1} and auxiliary variables $u_{i,1}^l = u_{i,0}^l = {x_{i,0}},\forall l \in \mathcal{Q}_i$, together with a gradient table $\left\{ {\nabla f_i^l\left( {u_{i,0}^l} \right)} \right\}_{l = 1}^{{q_i}}$.}
	\For{$k=0,1,\ldots,$}{
		\For{\rm{each reliable agent} $i \in \mathcal{R}$}{
			{\textbf{Exchange} information according to Step 3 in Algorithm \ref{Algo1}}.\\
			{\textbf{Select} uniformly a random sample with index $s_{i,k}$ from the set $\mathcal{Q}_i$ and evaluate the local stochastic gradient according to}
                \begin{equation*}\label{E3-3}
                {r_{i,k}^u} \!=\! \nabla f_i^{{s_{i,k}}}\left( {{x_{i,k}}} \right) - \nabla f_i^{{s_{i,k}}}\left( {u_{i,k}^{{s_{i,k}}}} \right)\! + \frac{1}{{{q_i}}}\sum\limits_{l = 1}^{{q_i}}\!\! {\nabla f_i^l\left( {u_{i,k}^l} \right)}.
                \end{equation*}\\
			{\textbf{Take} $u_{i,k + 1}^{s_{i,k}} = {x_{i,k}}$ and replace $\nabla f_i^{{s_{i,k}}}\left( {u_{i,k + 1}^{s_{i,k}}} \right) $ by $\nabla f_i^{{s_{i,k}}}\left( {{x_{i,k}}} \right)$ in the corresponding position of the gradient table, while keep the rest of the positions unchanged, i.e., $\nabla f_i^l\left( {u_{i,k + 1}^l} \right) = \nabla f_i^l\left( {u_{i,k}^l} \right)$, $ l \in {\mathcal{Q}_i}\backslash \left\{ {{s_{i,k}}} \right\}$;}\\
            {\textbf{Update} its current model according to Steps 7-8 in Algorithm \ref{Algo1}.}
		}
		\For{\rm{each Byzantine agent} $i \in \mathcal{B}$}{\textbf{Send} an arbitrary vector $ {z_{ij,k}}$ to its neighbor $j$, $j \in {\mathcal{N}_i}$.}
	}
	\KwOut{all decision variables ${x_{i,k}}$, $i \in \mathcal{V}$, until a prescribed criterion is satisfied.}
	\caption{\textit{Prox-DBRO-SAGA}.}
	\label{Algo2}
\end{algorithm}

\begin{algorithm}[!htp]
	\small
    \SetKwInput{KwInit}{Initialize}
	\SetKwBlock{Repeat}{Repeat for $ k = 0,1,2, \dots$}{}
	\SetKwBlock{End}{End for a required criterion is satisfied}{}
	\KwIn{a proper constant or decaying step-size $\alpha_k>0$.}
	\KwInit{the same parameters and starting points according to Algorithm \ref{Algo1} and an auxiliary variable ${w_{i,0}} = {x_{i,0}}$.}
	\For{$k=0,1,\ldots,$}{
		\For{\rm{each reliable agent} $i \in \mathcal{R}$}{
			{\textbf{Exchange} information according to Step 3 in Algorithm \ref{Algo1}}.\\
			{\textbf{Select} uniformly a random sample with index $s_{i,k}$ from the set $\mathcal{Q}_i$ and evaluate the local stochastic gradient according to}
            \begin{equation*}\label{E3-4}
            {r_{i,k}^w}= \nabla f_i^{{s_{i,k}}}\left( {{x_{i,k}}} \right) -  \nabla f_i^{{s_{i,k}}}\left( {{w_{i,k}}} \right) + \frac{1}{{{q_i}}}\sum\limits_{l = 1}^{{q_i}} {\nabla f_i^l\left( {{w_{i,k}}} \right)}.
            \end{equation*}\\
			{\textbf{Take} $w_{i,k+1} = {x_{i,k}}$ with a heterogenous triggering probability $p_i$ and keep $w_{i,k+1} = {w_{i,k}}$ with the probability $1-p_i$;}\\
            {\textbf{Update} its current model according to Steps 7-8 in Algorithm \ref{Algo1}.}
		}
		\For{\rm{each Byzantine agent} $i \in \mathcal{B}$}{\textbf{Send} an arbitrary vector $ {z_{ij,k}}$ to its neighbor $j$, $j \in {\mathcal{N}_i}$.}
	}
	\KwOut{all decision variables ${x_{i,k}}$, $i \in \mathcal{V}$, until a prescribed criterion is satisfied.}
	\caption{\textit{Prox-DBRO-LSVRG}.}
	\label{Algo3}
\end{algorithm}
\begin{remark}\label{R1}
All steps in Algorithms \ref{Algo1}-\ref{Algo3} are executed in parallel among all reliable agents since they are honest and hence comply with these update rules. It is also  worthwhile to mention that the expected cost in evaluating the local stochastic gradient under \textit{Prox-DBRO-LSVRG} is at least double that of \textit{Prox-DBRO-SAGA} at every iteration. This computational advantage of \textit{Prox-DBRO-SAGA} is at the expense of a higher (total) storage cost of ${\mathcal{O}}\left( {nQ} \right)$ (with $Q = \sum\nolimits_{i = 1}^{m}{{q_i}}$) than ${\mathcal{O}}\left( mn \right)$ of \textit{Prox-DBRO-LSVRG}, which is prohibitively expansive when the number ($q_i$) of local training samples is massive. Therefore, adopting either \textit{Prox-DBRO-SAGA} or \textit{Prox-DBRO-LSVRG} in practice involves a trade-off between per-iteration computational cost and storage. Users may improve and implement \textit{Prox-DBRO-VR} via incorporating other categories of VR techniques \cite{Gorbunov2019a} based on their customized needs.
\end{remark}

\section{Convergence Analysis}\label{ConAna}
To streamline notation, we denote ${{\mathcal{F}_k}}$ as the filter of the history w.r.t. the dynamical system generated by the sequence ${\left\{ {{s_{i,k}}} \right\}_{i \in \mathcal{R},k \ge 0}}$, and the conditional expectation $\mathbb{E}\left[ {{s_k}|{\mathcal{F}_k}} \right]$ is shortly denoted by ${\mathbb{E}_k}\left[ \cdot \right]$ in the sequel analysis. Let ${x_k}: = {\mathbf{col}}{\left\{ {{x_{i,k}}} \right\}_{i \in \mathcal{R}}} \in {\mathbb{R}^{\left| \mathcal{R} \right|n}}$ and ${r_k}: = {\mathbf{col}}{\left\{ {{r_{i,k}}} \right\}_{i \in \mathcal{R}}} \in {\mathbb{R}^{\left| \mathcal{R} \right|n}}$.
To facilitate the subsequent analysis, we define
\begin{equation*}
\begin{aligned}
&\nabla F\left( {{x_k}} \right): = {\mathbf{col}}{\left\{ {\nabla {f_i}\left( {{x_{i,k}}} \right)} \right\}_{i \in \mathcal{R}}},\\
&{\chi _i}\left( {{x_{i,k}}} \right): = \phi \sum\limits_{j \in {\mathcal{R}_i}} {{{\left\| {{x_{i,k}} - {x_{j,k}}} \right\|}_a}} ,\forall i \in \mathcal{R}, \\
&{\partial _x}\chi \left( {{x_k}} \right): = {\mathbf{col}}{\left\{ {{\partial _{{x_i}}}{\chi _i}\left( {{x_{i,k}}} \right)} \right\}_{i \in \mathcal{R}}},\\
&{\delta _i}\left( {{x_{i,k}}} \right): = \phi \sum\limits_{j \in {\mathcal{B}_i}} {{{\left\| {{x_{i,k}} - {x_{j,k}}} \right\|}_a}} ,\forall i \in \mathcal{R},\\
&{\partial _x}\delta \left( {{x_k}} \right): = {\mathbf{col}}{\left\{ {{\partial _{{x_i}}}{\delta _i}\left( {{x_{i,k}}} \right)} \right\}_{i \in \mathcal{R}}},\\
\end{aligned}
\end{equation*}
such that \textit{Prox-DBRO-VR} can be concisely formulated as:
\begin{subequations}\label{E3-5}
\begin{align}
\label{E3-5-1}&{\bar x_{k}} =  {x_k} - {\alpha _k}\left( {r_k }+ {\partial _x} \chi \left( {{x_k}} \right) {  + {\partial _x} \delta \left( {{x_k}} \right)} \right),\\
\label{E3-5-2}&{x_{k + 1}} =  {\mathbf{prox}}{_{{\alpha _k}, G}}\left\{ {{{\bar x}_k}} \right\}.
\end{align}
\end{subequations}
\subsection{Auxiliary Results}\label{AuxiRe}
Inspired by the unified analysis framework for centralized stochastic gradient descent methods in \cite{Gorbunov2019a}, we introduce the following two lemmas. To begin with, we define respectively two sequences for \textit{\textit{Prox-DBRO-SAGA}} and \textit{\textit{Prox-DBRO-LSVRG}} in the following. For \textit{\textit{Prox-DBRO-SAGA}}, we define
\begin{equation*}
t_{i,k}^u: = \frac{1}{{{q_i}}}\sum\limits_{l = 1}^{{q_i}} { {f_i^l\left( {u_{i,k}^l} \right) - f_i^l\left( {{{\tilde x}^*}} \right) - \nabla f_i^l{{\left( {{{\tilde x}^*}} \right)}^ \top }\left( {u_{i,k}^l - {{\tilde x}^*}} \right)} }.
\end{equation*}
For \textit{\textit{Prox-DBRO-LSVRG}}, we define
\begin{equation*}
t_{i,k}^w: = \frac{1}{{{q_i}}}\sum\limits_{l = 1}^{{q_i}} {{f_i^l\left( {{w_{i,k}}} \right) - f_i^l\left( {{{\tilde x}^*}} \right) - \nabla f_i^l{{\left( {{{\tilde x}^*}} \right)}^ \top }\left( {{w_{i,k}} - {{\tilde x}^*}} \right)} }.
\end{equation*}
Note that both sequences ${\left\{ {t_{i,k}^u} \right\}_{i \in \mathcal{R},k \ge 0}}$ and ${\left\{ {t_{i,k}^w} \right\}_{i \in \mathcal{R},k \ge 0}}$ are non-negative according to the convexity of the local component function ${f_i^l}$, $l \in \mathcal{Q}_i$. For the sequel analysis, we define respectively the gradient-learning quantities $t_k^u: = \sum\nolimits_{i \in \mathcal{R}} {t_{i,k}^u} $ and $t_k^w: = \sum\nolimits_{i \in \mathcal{R}} {t_{i,k}^w} $ for \textit{Prox-DBRO-SAGA} and \textit{Prox-DBRO-LSVRG}, the largest and smallest number of local samples ${q_{\min }} : = {\min _{i \in \mathcal{R}}}{q_i}$ and ${q_{\max }}: = {\max _{i \in \mathcal{R}}}{q_i}$, the minimum and maximum triggering probabilities ${p_{\min }} : = {\min _{i \in \mathcal{R}}}{p_i}$ and ${p_{\max }}: = {\max _{i \in \mathcal{R}}}{p_i}$, while ${\kappa _q}: = {q_{\max }}/{q_{\min }} \ge 1$.

\begin{lemma}(Gradient-learning sequences)\label{L1}
Since $t_k^u$ and $t_k^w$ are non-negative according to the convexity of the local component function ${f_i^l}$, $l \in \mathcal{Q}_i$ under Assumption \ref{Assu1}, $\forall k \ge 0$, we have for \textit{\textit{Prox-DBRO-SAGA}},
\begin{equation}\label{E4-1}
{\mathbb{E}_k}\left[ {{t_{k + 1}^u}} \right] \le \left( {1 - \frac{1}{{{q_{\max }}}}} \right)t_k^u + \frac{{D_F}\left( {{x_k},{x^*}} \right)}{{{q_{\min }} }},
\end{equation}
and for \textit{\textit{Prox-DBRO-LSVRG}},
\begin{equation}\label{E4-2}
\begin{aligned}
{\mathbb{E}_k}\left[ {t_{k + 1}^w} \right] \le  \left( {1 - {p_{\min }} } \right)t_k^w + {p_{\max }}{D_F}\left( {{x_k},{x^*}} \right),
\end{aligned}
\end{equation}
where ${D_F}\left( {{x_k},{x^*}} \right): = F\left( {{x_k}} \right) - F\left( {{x^*}} \right) - \nabla F{\left( {{x^*}} \right)^ \top }\left( {{x_k} - {x^*}} \right)$ is known as the Bregman divergence with respect to the convex cost function $F$.
\begin{proof}
See Appendix \ref{Appen1}.
\end{proof}
\end{lemma}
We next seek the upper bound of the distance between the local stochastic gradient estimator $r_k$ and gradient $\nabla F\left( {{x^*}} \right)$ at the optimal solution for both \textit{\textit{Prox-DBRO-SAGA}} and \textit{\textit{Prox-DBRO-LSVRG}}.
\begin{lemma}(Gradient-learning errors)\label{L2}
Suppose that Assumptions \ref{Assu1}-\ref{Assu2} hold. For $ k \ge 0$, we have for \textit{\textit{Prox-DBRO-SAGA}},
\begin{equation}\label{E4-3}
{\mathbb{E}_k}\left[ {{{\left\| {{r_k^u} - \nabla F\left( {{x^*}} \right)} \right\|}_2^2}} \right] \le 4Lt_k^u + 2\left( {2L - \mu } \right){D_F}\left( {{x_k},{x^*}} \right),
\end{equation}
and for \textit{Prox-DBRO-LSVRG},
\begin{equation}\label{E4-4}
{\mathbb{E}_k}\left[ {{{\left\| {{r_k^w} - \nabla F\left( {{x^*}} \right)} \right\|}_2^2}} \right] \le 4Lt_k^w + 2\left( {2L - \mu } \right){D_F}\left( {{x_k},{x^*}} \right),
\end{equation}
where $r_k^u: = {\mathbf{col}}{\left\{ {r_{i,k}^u} \right\}_{i \in \mathcal{R}}} \in {\mathbb{R}^{\left| \mathcal{R} \right|n}}$ and $r_k^w: = {\mathbf{col}}{\left\{ {r_{i,k}^w} \right\}_{i \in \mathcal{R}}} \in {\mathbb{R}^{\left| \mathcal{R} \right|n}}$.
\begin{proof}
See Appendix \ref{Appen2}.
\end{proof}
\end{lemma}
The following proposition is an important result for the analysis of arbitrary norm approximation.
\begin{proposition}\label{P1}
Consider two constants $a_1 \ge 1$ and $a_2$, such that $1/{a_1} + 1/{a_2} = 1$. For an arbitrary vector $\tilde x \in {\mathbb{R}^n}$, we denote the subdifferential ${\partial}{\left\| {\tilde x} \right\|_{a_1}} = \left\{ {\tilde z \in {\mathbb{R}^n}:\left\langle {\tilde z,\tilde x} \right\rangle  = {{\left\| {\tilde x} \right\|}_{a_1}},{{\left\| {\tilde z} \right\|}_{a_2}} \le 1} \right\}$.
\begin{proof}
We refer interested readers to the supplementary document of \cite{Li2019k} for the proof of Proposition \ref{P1}.
\end{proof}
\end{proposition}

\begin{proposition}\label{P2}
Recalling the definition of ${\mathbf{prox}}_{\alpha ,g}\left\{ {x_i} \right\}$, we know that ${\left[{\mathbf{prox}}_{\alpha, G}\left\{ x \right\}\right]}_i = {\mathbf{prox}}_{\alpha ,g}\left\{ {x_i} \right\}$,  $\forall i \in \mathcal{R}$, and
\begin{equation}\label{E4-5}
{\left\| {{\mathbf{prox}}_{\alpha, G}\left\{ x \right\} - {\mathbf{prox}}_{\alpha, G}\left\{ y \right\}} \right\|_2} \le {\left\| {x - y} \right\|_2},
\end{equation}
where $x = {\mathbf{col}}{\left\{ {{x_i}} \right\}_{i \in \mathcal{R}}} \in {\mathbb{R}^{\left| \mathcal{R} \right|n}}$ and $y = {\mathbf{col}}{\left\{ {{y_i}} \right\}_{i \in \mathcal{R}}} \in {\mathbb{R}^{\left| \mathcal{R} \right|n}}$.
\begin{proof}
See Appendix \ref{Appen3}.
\end{proof}
\end{proposition}

\subsection{Main Results}\label{MainRe}
We next derive a feasible range for the penalty parameter to enable the equivalence between the decentralized consensus optimization problem (\ref{E2-3}) and norm-penalized approximation problem (\ref{E2-4}) as follows, which further guarantees the equivalence between the original optimization problem (\ref{E2-1}) and norm-penalized approximation problem (\ref{E2-4}). We define $\Pi  \in {\mathbb{R}^{\left| \mathcal{R} \right| \times \left| {{\mathcal{E}_\mathcal{R}}} \right|}}$ by the incidence matrix associated with ${\mathcal{G}_\mathcal{R}}$. Specifically, for any edge $e:=\left( {i,j} \right) \in {{\mathcal{E}_\mathcal{R}}}$ with $i<j$, the $\left( {i,e} \right)$-th and $\left( {j,e} \right)$-th entries of $\Pi$ are $1$ and $-1$, respectively.
\begin{theorem}(Resilient consensus condition)\label{T1}
Suppose that Assumptions \ref{Assu1} and \ref{Assu2} hold. Given $g'\left( {{{\tilde x}^*}} \right) \in {\partial}{g}\left( {{{\tilde x}^*}} \right)$, if the penalty parameter satisfy $\phi  \ge {\phi _{\min }}: = {\left| \mathcal{R} \right|^{\frac{3}{2}}}\sqrt {\left| {{\mathcal{E}_\mathcal{R}}} \right|} \mathop {\max }\nolimits_{i \in \mathcal{R}} {\left\| {\nabla {f_i}\left( {{{\tilde x}^*}} \right) }+ {g'\left( {{{\tilde x}^*}} \right)} \right\|_\infty }/{\lambda _{\min }}\left( \Pi  \right)$, the optimal solution to the original optimization problem (\ref{E2-1}) is equivalent to the globally optimal solution to norm-penalized approximation problem (\ref{E2-4}), i.e., ${x^*} = {1_{\left| \mathcal{R} \right|}} \otimes {{\tilde x}^*}$.
\begin{proof}
See Appendix \ref{Appen4}.
\end{proof}
\end{theorem}
\begin{remark}\label{R4-1}
Theorem \ref{T1} demonstrates that a selection of a sufficiently large penalty parameter guarantees the equivalence between the original optimization problem (\ref{E2-1}) and norm-penalized approximation problem (\ref{E2-4}). However, the sequel convergence results manifest that a larger $\phi$ causes a bigger convergence error. Therefore, the notion of a sufficiently large penalty parameter is tailored for theoretical results, and one can hand-tune this parameter to obtain better algorithm performances in practice.
\end{remark}

For simplicity, we fix the minimum and maximum heterogenous triggering probabilities as ${p_{\min }} = 1/{q_{\max }}$ and $p_{\max} = 1/ {q_{\min }}$, respectively. Hence, we define ${r_k}: = \left\{ \begin{gathered}
  r_k^u,{\text{for \textit{Prox-DBRO-SAGA}}} \hfill \\
  r_k^w,{\text{for \textit{Prox-DBRO-LSVRG}}} \hfill \\
\end{gathered}  \right.$ and ${t_k}: = \left\{ \begin{gathered}
  t_k^u,\text{for \textit{Prox-DBRO-SAGA}} \hfill \\
  t_k^w,\text{for \textit{Prox-DBRO-LSVRG}} \hfill \\
\end{gathered}  \right.$, such that the theoretical results for both \textit{Prox-DBRO-SAGA} (Algorithm \ref{Algo2}) and \textit{Prox-DBRO-LSVRG} (Algorithm \ref{Algo3}) can be unified in a general framework. Before deriving a linear convergence rate for Algorithms \ref{Algo2}-\ref{Algo3}, we first define the sequel parameters: $\gamma : = \mu L/\left( {\mu  + L} \right)$, $P_1: = {4n{\phi ^2}\sum\nolimits_{i \in \mathcal{R}} {\left( {2{{\left| {{\mathcal{R}_i}} \right|}^2} + {{\left| {{\mathcal{B}_i}} \right|}^2}} \right)} }$, $P_2: = n{\phi ^2}\sum\nolimits_{i \in {\mathcal R}} {{{\left| {{{\cal B}_i}} \right|}^2}} /\gamma $, and $E := 4{P_2}/\gamma $.
\begin{theorem}(Linear convergence).\label{T2}
Suppose that Assumptions \ref{Assu1}-\ref{Assu2} hold. Under the conditions of Theorem \ref{T1},  if the constant step-size meets $0 < {\alpha _k} \equiv \alpha  \le 1/\left( {{\kappa _q}\left( {32{{\left( {1 + {\kappa _f}} \right)}^2} + q_{\text{min}}} \right)\mu} \right)$, then the sequence ${\left\{ {{x_k}} \right\}_{k \ge 0}}$ generated by Algorithms \ref{Algo2}-\ref{Algo3}, converges linearly to an error ball around the optimal solution to the original optimization problem (\ref{E2-1}) at a linear rate of ${\left( {1 - \mathcal{O}\left( {\gamma \alpha } \right)} \right)^k}$, i.e.,
\begin{equation}\label{E4-6}
\begin{aligned}
&\mathbb{E}\left[ {\left\| {{x_k} - {1_{\left| \mathcal{R} \right|}} \otimes {{\tilde x}^*}} \right\|_2^2} \right]\\
\le & {\left( {1 - \frac{\gamma }{4}\alpha } \right)^k}{U_0} + 4\left( {\frac{{{P_1}}}{\gamma }\alpha  + E} \right)\left( {1 - {{\left( {1 - \frac{\gamma }{4}\alpha } \right)}^{k}}} \right),
\end{aligned}
\end{equation}
where ${U_0} = \left\| {{x_0} - {x^*}} \right\|_2^2 + {{q_{\min }}}\gamma \alpha {t_0}/\left( {{{q_{\max }}}L} \right)$, and the radius of the error ball is no more than $4\left( {{P_1}\alpha /\gamma  + E} \right)$.
\begin{proof}
See Appendix \ref{Appen5}.
\end{proof}
\end{theorem}


To proceed, we define $\theta  > 4/\gamma $, $\Xi : = \max \left\{ {\frac{{{\theta ^2}P_1}}{{\gamma \theta  - 1}},\left( {\xi  - \frac{\gamma }{4}\theta } \right) {\left\| {{x_0} - {x^*}} \right\|_2^2} + \frac{{{\theta ^2}}}{\xi }P_1 + \theta {P_2} - \xi E} \right\}$, and $\xi := {\kappa _q}\left( {64{{\left( {1 + {\kappa _f}} \right)}^2} + {q_{\min }}} \right)\mu \theta $.
\begin{theorem}(Sub-linear Convergence).\label{T3}
Suppose that Assumptions \ref{Assu1}-\ref{Assu2} hold. Under the condition of Theorem \ref{T1}, if the decaying step-size is chosen as ${\alpha _k} = \theta /\left( {k + \xi} \right)$, then the sequence ${\left\{ {{x_k}} \right\}_{k \ge 0}}$ generated by Algorithms \ref{Algo2}-\ref{Algo3} converges to an error ball around the optimal solution to the original optimization problem (\ref{E2-1}), at a sub-linear rate of $\mathcal{O}\left( {1/k } \right)$, i.e.,
\begin{equation}\label{E4-7}
\begin{aligned}
\mathbb{E}\left[ {\left\| {{x_k} - {1_{\left| \mathcal{R} \right|}} \otimes {{\tilde x}^*}} \right\|_2^2} \right] \le \frac{\Xi }{{k + \xi}} + E,\forall k \ge 0,
\end{aligned}
\end{equation}
where the radius of the error ball is $E$.
\begin{proof}
See Appendix \ref{Appen6}.
\end{proof}
\end{theorem}
\begin{remark}\label{R4-1}
The convergence results established in Theorems \ref{T2}-\ref{T3} assert that the proposed algorithms achieve linear convergence at the expense of a larger
convergence error than that of in the sub-linear convergence case. We note that a smaller constant step-size in the linear convergence case may simultaneously lead to a smaller convergence error and a slower convergence rate according to Theorem \ref{T2}. It is also clear from Theorem \ref{T3} that the convergence error of \textit{Prox-DBRO-SAGA} and \textit{Prox-DBRO-LSVRG} for the sub-linear convergence case is determined by the number of Byzantine agents. That is to say, the sub-linear exact convergence of \textit{Prox-DBRO-SAGA} and  \textit{Prox-DBRO-LSVRG} can be recovered, when the number of Byzantine agents equals to zero. Therefore, the theoretical results derived in Theorems \ref{T2}-\ref{T3} demonstrate a trade-off between the convergence error and convergence rate, which has been revealed by the theoretical results.
\end{remark}

\begin{remark}\label{R4-2}
This paper does not make any assumption or restrictions on the number/proportion of Byzantine agents in the network and only assumes a connected network among all reliable agents (see Assumption \ref{Assu2}). However, this does not imply that the number of Byzantine agents can be unbounded since an increase in the number of Byzantine agents lead to larger convergence errors in both cases and an unbounded number of Byzantine agents causes eventually divergence of the proposed algorithms according to Eqs. (\ref{E4-6}) and (\ref{E4-7}) in Theorems 2-3.
According to Theorems \ref{T2}-\ref{T3}, the resilience of Algorithms \ref{Algo2}-\ref{Algo3} is characterized by the consensus and controllable convergence errors of all reliable agents.
\end{remark}
The sequel corollary verifies the statement that there is no restriction on the number or proportion of Byzantine agents under appropriate assumptions.
\begin{corollary}(Guaranteed convergence for a high proportion of Byzantine agents)\label{C1}
Suppose that Assumptions \ref{Assu1}-\ref{Assu2} hold. Under the conditions of Theorem \ref{T2} (resp., Theorem \ref{T3}), if we set the number of Byzantine agents by $\left| {\cal B} \right| = {r_a}\left| {\cal R} \right|$ with any bounded constant $0< r_a < \infty$ that is only need to be properly chosen to avoid floating point issues, then the sequence ${\left\{ {{x_k}} \right\}_{k \ge 0}}$ generated by Algorithms \ref{Algo2}-\ref{Algo3} with the constant (resp., decaying) step-size, still converges to an error ball around the optimal solution to the original problem (\ref{E2-1}) at a linear (resp., sub-linear) rate of ${\left( {1 - \mathcal{O}\left( {\gamma \alpha } \right)} \right)^k}$ (resp., $\mathcal{O}\left( {1/k } \right)$), $\forall k \ge 1$, where the radius of the error ball is bounded.
\begin{proof}
We begin by showing that the convergence error is bounded for any $0< r_a < \infty$ in the linear convergence case. By relaxing the upper bound in the right-hand-side (RHS) of  (\ref{E4-6}), we have
\begin{equation}\label{CE-1}
\begin{aligned}
&\mathbb{E}\left[ {\left\| {{x_k} - {1_{\left| {\cal R} \right|}} \otimes {{\tilde x}^*}} \right\|_2^2} \right] \\
\le & {\left( {1 - \frac{\gamma }{4}\alpha } \right)^k}{U_0} + {4n\left( {4\alpha \left( {\frac{{4 + r_a^2}}{\gamma }} \right) + {{\left( {\frac{{{r_a}}}{\gamma }} \right)}^2}} \right){{\phi }^2}{{\left| {\cal R} \right|}^3}}.
\end{aligned}
\end{equation}
where the first term in the RHS of (\ref{CE-1}) decays linearly over the iteration $k$, $\forall k \ge 1$, as $0 < 1 - \gamma \alpha /4 < 1$ according to Theorem \ref{T2} and $U_0$ is a constant owing to the initialization rule of the proposed algorithms. Moreover, since both $n$, $\gamma $, $\phi $, and the step-size $\alpha$ are (bounded) constants once the objective function and network are determined, we know that the second term, i.e., the error term, in the RHS of (\ref{CE-1}) can be bounded by a finite value when $0< r_a < \infty$. Therefore, the radius of the convergence error ball is bounded. We omit the proof for the sub-linear case since it follows a same technical line as in the linear case.
\end{proof}
\end{corollary}
\begin{remark}\label{R4-2}
Corollary \ref{C1} asserts that the radius of the error ball is bounded without any restriction on the number or proportion of Byzantine agents under some appropriate conditions. A specific example is when the proportion of Byzantine agents accounts for 50\%, i.e., $r_a = 1$, the radius of the convergence error ball is $ 4n\left( {20\gamma \alpha  + 1} \right){{\phi }^2}{\left| {\cal R} \right|^3}/{\gamma ^2}$ according to (\ref{CE-1}), which is bounded since all the parameters are constants once the objective function and network are determined. However, this does not imply that the number of Byzantine agents can be arbitrarily large since an increase in the number of Byzantine agents lead to larger convergence errors in both cases and an unbounded number, i.e., ${r_a} \to \infty $, of Byzantine agents causes eventually divergence of the proposed algorithms according to Theorems \ref{T2}-\ref{T3}
\end{remark}
\section{Numerical Experiments}\label{ExpRe}
In this section, we perform a case study on decentralized soft-max regression with sparsity to verify the theoretical results and show the convergence performance of the proposed algorithms, where four kinds of Byzantine attacks (zero-sum attacks, Gaussian attacks, same-value attacks, and sign-flipping attacks) are considered. The communication networks are randomly generated by the Erdős-Rényi method, where Byzantine agents are also selected in a random way.
Most existing literature adopts only testing accuracy and the consensus error to validate the convergence performance of their proposed algorithms. However, neither higher testing accuracy nor a smaller consensus error can comprehensively reflect the convergence of the tested algorithms. This is because these two metrics fail to precisely measure the iterative distance between the optimized function value and the optimal value of optimization problems, which, however, serves as a primary goal of theoretical analysis. Therefore, there is a gap between the theoretical result and practical performance. To bridge the gap, we introduce the (averaged) optimal gap in the form of function values, i.e., $\left( {1/\left| \mathcal{R} \right|} \right)\sum\nolimits_{i \in \mathcal{R}} {\left( {{f_i}\left( {{x_{i,k}}} \right) + g\left( {{x_{i,k}}} \right) - \left( {{f_i}\left( {{{\tilde x}^*}} \right) + g\left( {{{\tilde x}^*}} \right)} \right)} \right)} $, as the third metric, which could precisely captures transient behaviors (convergence or divergence) of algorithms when training a machine- or deep-learning model. A network of $m$ agents consisting of ${\left| \mathcal{R} \right|}$ reliable agents and $\left| \mathcal{B} \right|$  Byzantine agents, minimize a regularized soft-max regression problem for a multi-class classification task via specifying the problem formulation (\ref{E2-1}) as
${f_i}\left( {\tilde x} \right): =  - \left( {1/{q_i}} \right)\sum\nolimits_{j = 1}^{{q_i}} {\sum\nolimits_{l = 0}^{\tilde C - 1} {{\mathcal{I}_{\left( {{\tilde b_{ij}} = l} \right)}}} } \ln \left( {{e^{\left[ {\tilde x} \right]_l^ \top {\tilde c_{ij}}}}/\sum\nolimits_{t = 0}^{\tilde C - 1} {{e^{\left[ {\tilde x} \right]_t^ \top {\tilde c_{ij}}}}} } \right) + \left( {{\beta _1}/2} \right)\left\| {\tilde x} \right\|_2^2$ and ${g}\left( {\tilde x} \right): = {\beta _2}{\left\| {\tilde x} \right\|_1}$, where $\tilde x \in {\mathbb{R}^{n}}$ with $n = \tilde C \tilde n$ is the model parameter, ${\tilde C}$ represents the number of sample classes, ${\left[ {\tilde x} \right]_j}$ denotes a vector that contains $\left(j \tilde n \right)$-th to $\left( \left( {j + 1} \right)\tilde n - 1\right)$-th elements of $\tilde x$, ${{{\tilde b_{ij}}}}$ and ${{\tilde c_{ij}}}$ are the $j$-th label and image allocated to agent $i$, respectively, ${{\mathcal{I}_{\left( {{{\tilde b_{ij}}} = l} \right)}}}$ is the indicator function with ${\mathcal{I}_{\left( {{{\tilde b_{ij}}} = l} \right)}} = 1$ if ${{{\tilde b_{ij}}} = l}$ and ${\mathcal{I}_{\left( {{{\tilde b_{ij}}} = l} \right)}} = 0$ otherwise, ${\beta _1}$ and ${\beta _2}$ are positive parameters of regularized terms for avoiding over-fitting and obtaining a sparse solution, respectively.
We denote the total number of training samples by $N$ and the regularized parameters are set as ${\beta _1} = {\beta _2} = 1/N$ in the following numerical experiments.
\begin{figure}[!htp]
  \centering
  \includegraphics[width=1.3in,height=1.3in]{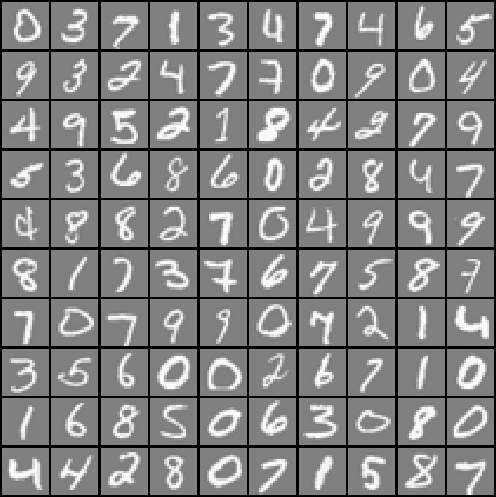}
  \caption{Random samples selected from the MNIST data set}\label{Fig-0}
\end{figure}
Since the algorithmic framework \textit{BRIDGE} \cite{Fang2022} and decentralized algorithm (denoted by \textit{Peng}) \cite{Peng2021} are only available to handling a class of smooth single-objective optimization problems, we equip them with the proximal-gradient mapping method \cite{Notarnicola2017a, Li2019m, Xu2021, Alghunaim2021, Li2021} to obtain \textit{Prox-BRIDGE-T}, \textit{Prox-BRIDGE-M}, \textit{Prox-BRIDGE-K}, and \textit{Prox-Peng} for the nonsmooth composite finite-sum optimization problem, which is also applied to \textit{GeoMed} \cite{Chen2018} to get \textit{Prox-GeoMed}. The initial state of decision variables of all tested algorithms are the same and generated from a standard normal distribution. Note that the parameters of all tested algorithms are optimized manually to obtain their best performance, and the parameters associated with the problem model keep the same to ensure fairness.
A total number of $Q=60000$ training samples from the MNIST \cite{LeCun2010} data set are evenly allocated to each agent (including both reliable agents and Byzantine agents in the network) to train the discriminator, while the rest $10000$ samples are used for testing. Fig. \ref{Fig-0} presents 100 samples randomly selected from the data set. Recall the theoretical results regarding the decaying and constant step-sizes and penalty parameter such that the experimental setting gives the following feasible selection ranges: $\alpha  \in \left( {0, 0.1385} \right]$, ${\alpha _k} \in \left( {0, 1/\left( {k + 14.2733} \right)} \right]$, and $\phi  \ge 0.0003$.
\begin{figure*}[!htp]
\begin{center}
\subfloat[A 30-agent network containing $5$ Byzantine agents.]{\includegraphics[width=1.7in,height=1.2in]{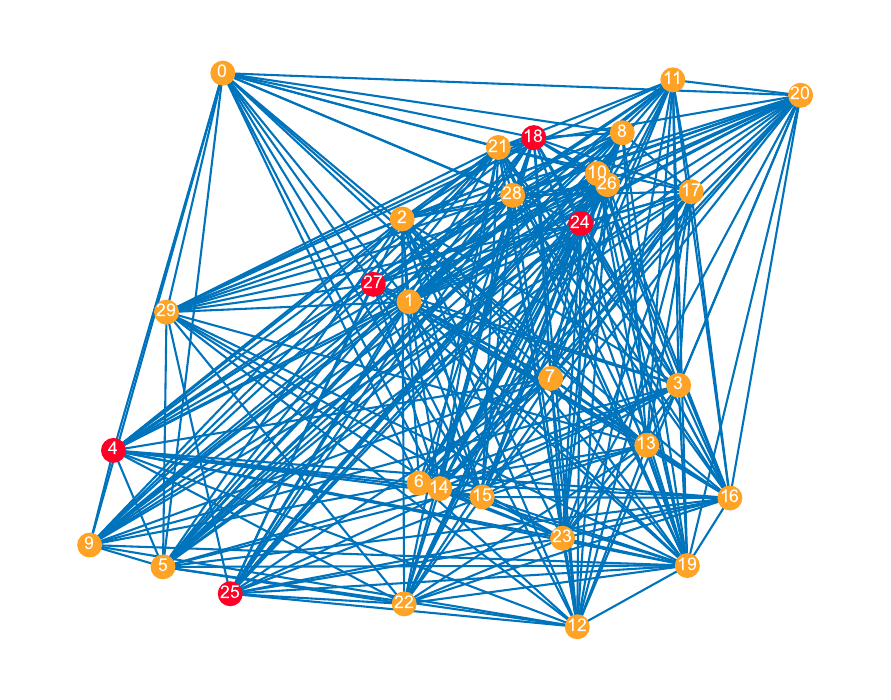}\label{Fig-1-1}} \hfill
\subfloat[Optimal gap over epochs.]{\includegraphics[width=1.7in,height=1.2in]{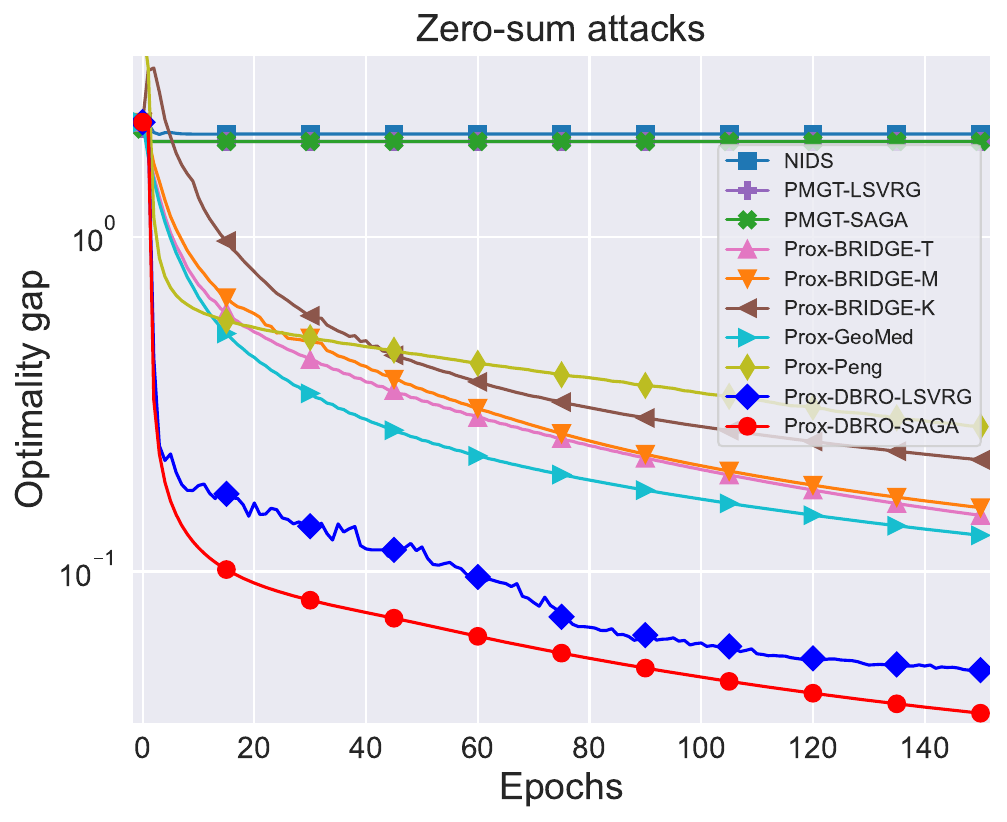}\label{Fig-1-2}}\hfill
\subfloat[Testing accuracy over epochs.]{\includegraphics[width=1.7in,height=1.2in]{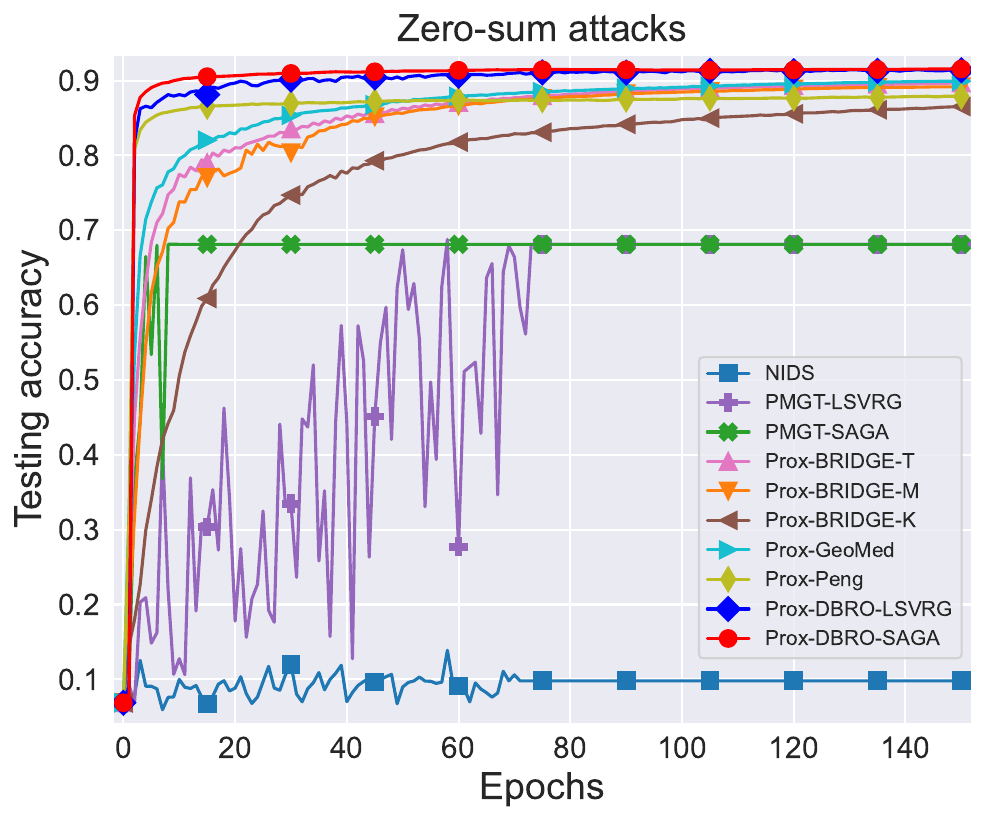}\label{Fig-1-3}}\hfill
\subfloat[Consensus error over epochs.]{\includegraphics[width=1.7in,height=1.2in]{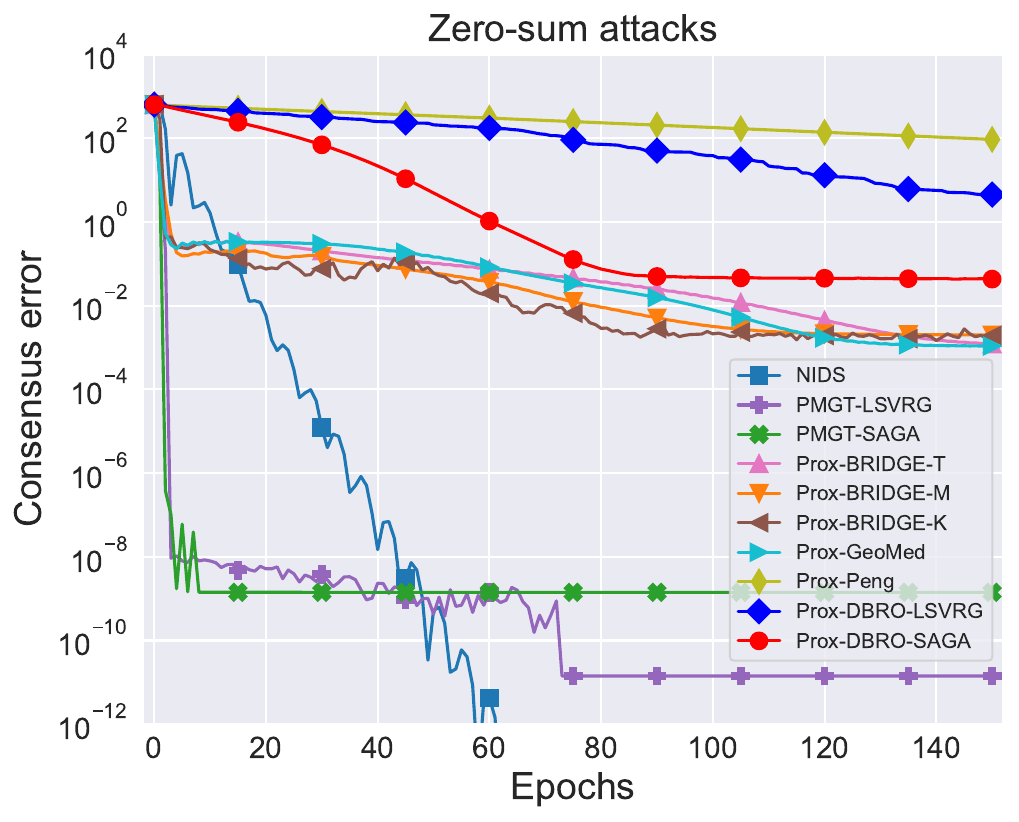}\label{Fig-1-4}}\hfill
\end{center}
\caption{Performance of all tested algorithms under zero-sum attacks.}
\label{Fig-1}
\end{figure*}
\begin{table*}[!htp]
\centering
\caption{Parameter settings and algorithm performance at 150 epochs under zero-sum attacks.}
\scalebox{0.70}{
\begin{tabular}{ccccccccccc}
\toprule
\diagbox{PS\&PI}{Algorithms} & \textit{NIDS} & \textit{PMGT-LSVRG} & \textit{PMGT-SAGA}  & \textit{Prox-BRIDGE-T} & \textit{Prox-BRIDGE-M} & \textit{Prox-BRIDGE-K} & \textit{Prox-GeoMed} & \textit{Prox-Peng} & {\bf{\textit{Prox-DBRO-LSVRG}}} & {\bf{\textit{Prox-DBRO-SAGA}}}\\
\midrule
Step-size & [0.01, 0.015] & 0.001 & 0.01  & 0.35  & 0.3 & 0.35 & 0.35 & 0.5 & 0.05 & 0.005  \\
Triggering probability & N/A & $m/Q$ & N/A  & N/A & N/A & N/A & N/A &N/A & $\left[ {m/Q/2,m/Q} \right]$ & N/A  \\
Penalty parameter & N/A & N/A & N/A & N/A & N/A & N/A & N/A & 0.1 & [0.2, 0.25] & [0.2, 0.25] \\
Consensus error & \bf{0} & 1.4163e-11 & 1.4159e-09 & 1.2492e-03 & 1.2116e-03 & 1.9340e-03 & 3.8010e-03 & 92.9523 & 4.4680 & 4.2898e-02   \\
Testing accuracy & 0.098 & 0.6812 & 0.6812 & 0.8965 & 0.8918 & 0.8653 & 0.8999 & 0.8789 & 0.9137 & \bf{0.9155} \\
Optimal gap & 2.0245 & 1.9230 & 1.9228 & 1.4940e-01 & 1.5751e-01 & 2.1701e-01 & 1.2435e-01 & 2.7077e-01 & 5.0790e-02 & \bf{3.7756e-02}   \\
\hline
\multicolumn{11}{c}{\makecell[c]{PS\&PI is the abbreviation of parameter settings and performance metrics.}}\\
\bottomrule
\end{tabular}}
\label{Tab-1}       
\end{table*}\newline
{\bf{Zero-sum attacks:}} As depicted in Fig. \ref{Fig-1}: \subref{Fig-1-1}, an $m=30$ network consists of ${\left| \mathcal{R} \right|}= 25$ reliable agents (yellow nodes) and $\left| \mathcal{B} \right| = 5$  Byzantine agents (red nodes), where each Byzantine agent $b$, $b \in \mathcal{B}$, sends a well-designed malicious message ${z_{ib,k}} =  - {\sum\nolimits_{j \in {\mathcal{R}_i}} {{w_{ij}}{x_{j,k}}} } /\left| {{\mathcal{B}_i}} \right|/{w_{bi}}$ to its reliable neighbor $i$, $i \in {\mathcal{R}_b}$, to drive the states of the reliable agent ${x_{i,k}} = 0_n$ at each iteration. For \textit{NIDS} \cite{Li2019m}, the algorithm parameter is fixed as $c = 1/\left( {2{{\max }_{i \in \mathcal{R}}}\left\{ {{\alpha _i}} \right\}} \right)$ for $\tilde W = {I_m} - c{D_\alpha }\left( {{I_m} - W} \right)$, where $\tilde W$, $W$, and $D_\alpha$ are the modified mixing matrix, weight matrix, and uncoordinated step-size, respectively. This means if $c$ is sufficiently small, \textit{NIDS} runs without any communication happening among all agents (both reliable agents and Byzantine agents) in the network. For \textit{PMGT-SAGA/PMGT-LSVRG} \cite{Ye2020}, we hand-tune the parameter associated with multi-step communications to obtain the best performance.
It is clear from Fig. \ref{Fig-1}: \subref{Fig-1-2}-\subref{Fig-1-4} and Table \ref{Tab-1} that the proposed algorithms achieve a smaller optimal gap and higher testing accuracy than the other tested algorithms in a same amount of computational costs (epochs). This demonstrates that the proposed algorithms approximate faster to the optimal solution than the other tested algorithms. Notably, zero-sum attacks launched by Byzantine agents aim to drive the states of all reliable agents to zero at each iteration. Therefore, a much smaller consensus error of \textit{NIDS} and \textit{PMGT-SAGA/PMGT-LSVRG} than the other Byzantine-resilient decentralized algorithms, may indicate that they are less resilient or more susceptible to the zero-sum attacks. This can be testified by their bigger optimal gaps and lower testing accuracy shown in Fig. \ref{Fig-1}: \subref{Fig-1-2}-\subref{Fig-1-3} and Table \ref{Tab-1}.
\begin{figure*}[!htp]
\begin{center}
\subfloat[A 40-agent network containing $8$ Byzantine agents.]{\includegraphics[width=1.7in,height=1.2in]{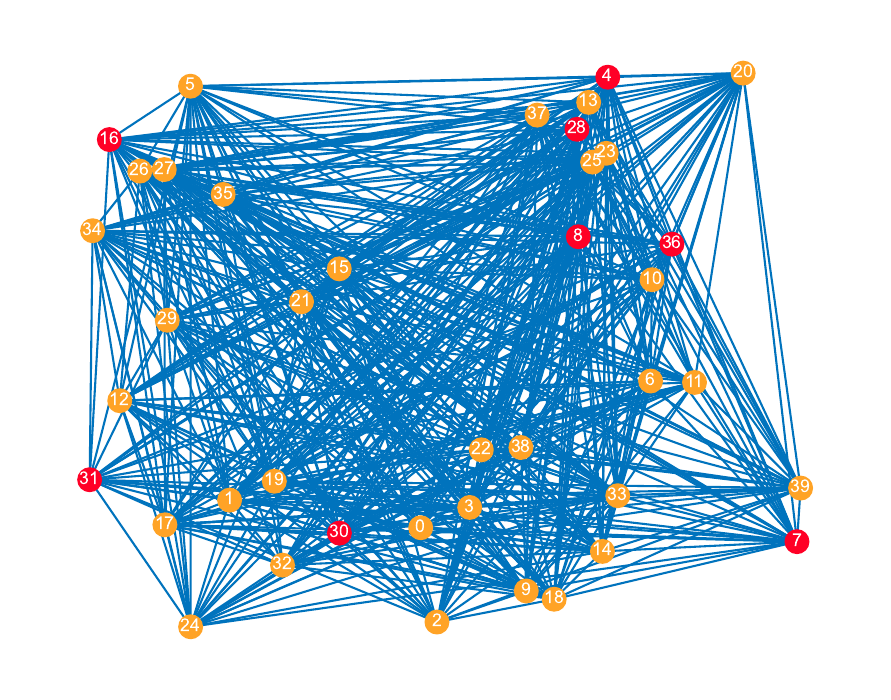}\label{Fig-2-1}} \hfill
\subfloat[Optimal gap over epochs.]{\includegraphics[width=1.7in,height=1.2in]{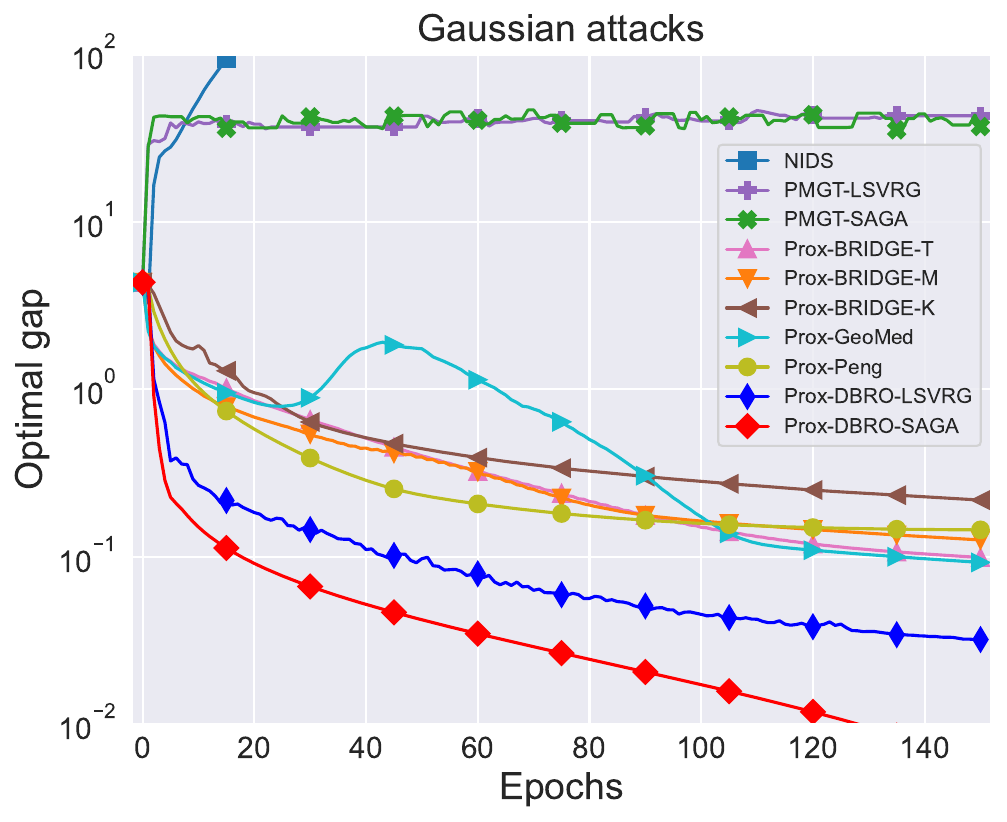}\label{Fig-2-2}} \hfill
\subfloat[Testing accuracy over epochs.]{\includegraphics[width=1.7in,height=1.2in]{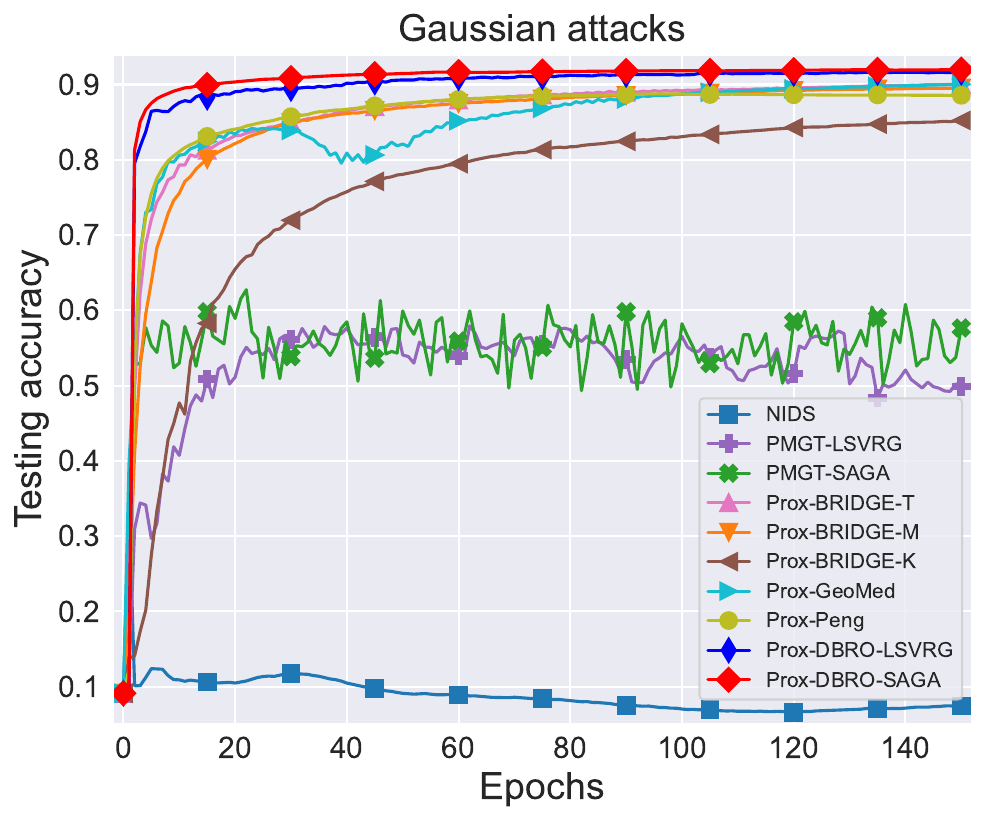}\label{Fig-2-3}} \hfill
\subfloat[Consensus error over epochs.]{\includegraphics[width=1.7in,height=1.2in]{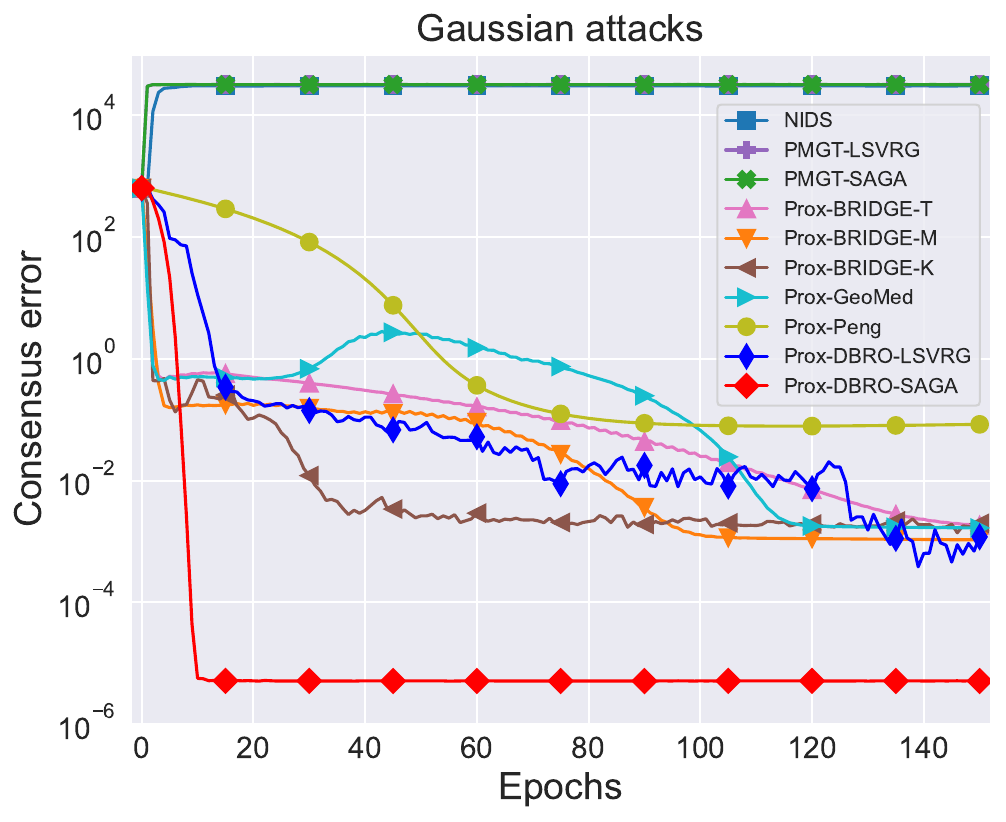}\label{Fig-2-4}} \hfill
\end{center}
\caption{Performance comparison among all tested algorithms under Gaussian attacks.}
\label{Fig-2}
\end{figure*}\newline
{\bf{Gaussian attacks:}} It is shown in Fig. \ref{Fig-2}: \subref{Fig-2-1} that an $m=40$ network consists of ${\left| \mathcal{R} \right|}= 32$ reliable agents (yellow nodes) and $\left| \mathcal{B} \right| = 8$  Byzantine agents (red nodes), where each Byzantine agent $b$, $b \in \mathcal{B}$, sends a message subject
to a Gaussian distribution with mean $ \sum\nolimits_{j \in {\mathcal{R}_i}} {{w_{ij}}{x_{j,k}}} /\sum\nolimits_{j \in {\mathcal{R}_i}} {{w_{ij}}} $ and standard
\begin{table*}[!htp]
\centering
\caption{Parameter settings and algorithm performance at 150 epochs under Gaussian attacks.}
\scalebox{0.70}{
\begin{tabular}{ccccccccccc}
\toprule
\diagbox{PS\&PI}{Algorithms} & \textit{NIDS} & \textit{PMGT-LSVRG} & \textit{PMGT-SAGA} & \textit{Prox-BRIDGE-T} & \textit{Prox-BRIDGE-M} & \textit{Prox-BRIDGE-K} & \textit{Prox-GeoMed} & \textit{Prox-Peng} & {\bf{\textit{Prox-DBRO-LSVRG}}} & {\bf{\textit{Prox-DBRO-SAGA}}}\\
\midrule
Step-size & [0.3, 0.35] & 0.3 & 0.3 & 0.4 & 0.3 & 0.3 & 0.4 & 0.05 & 0.0015 & 0.0025 \\
Triggering probability & N/A & $m/Q$ & N/A & N/A & N/A & N/A & N/A & N/A & $\left[ {m/Q/4,m/Q/2} \right]$ & N/A  \\
Penalty parameter & N/A & N/A & N/A & N/A & N/A & N/A & N/A & 0.001 & [0.001, 0.0015] & [0.001, 0.0015] \\
Consensus error & 3.0396e+04 & 3.2041e+04 & 3.1770e+04 & 1.8532e-03 & 1.0682e-03 & 1.9835e-03 & 1.6716e-03 & 8.3993e-02 & 1.1938e-03 & \bf{5.1252e-06} \\
Testing accuracy & 0.0742 & 0.4985 & 0.5765 & 0.9001 & 0.8952 & 0.8522 & 0.9002 & 0.8856 & 0.9165 & \bf{0.9197} \\
Optimal gap & Inf & 4.2709e+01 & 3.7046e+01 & 1.2136e-01 & 1.4870e-01 & 2.4074e-01 & 1.1507e-01 & 1.4480e-01 & 3.2053e-02 & \bf{5.7316e-03} \\
\hline
\multicolumn{11}{c}{\makecell[c]{PS\&PI is the abbreviation of parameter settings and performance metrics.}}\\
\bottomrule
\end{tabular}}
\label{Tab-2}       
\end{table*}
deviation 30, to its reliable neighbor $i$, $i \in {\mathcal{R}_b}$ at each iteration. This attack serves as a Gaussian noise, which can easily inflict
fluctuation on the state of reliable agents and deviate the states from their true values. Even though the testing accuracy index can still fluctuate around 0.6, we can see from Fig. \ref{Fig-2}: \subref{Fig-2-2} and Table \ref{Tab-2} that \textit{NIDS} and \textit{PMGT-SAGA/PMGT-LSVRG} show divergence from the optimal solution under Gaussian attacks. It is shown by Figs. \ref{Fig-2}: \subref{Fig-2-2}-\subref{Fig-2-4} that the proposed algorithms can still achieve a smaller optimal gap and higher testing accuracy in the same epochs, alternatively faster convergence, than the other tested algorithms. Moreover, one can clearly see from Table \ref{Tab-2} that \textit{Prox-DBRO-SAGA} takes the superiority on all three performance metrics (optimal gap, testing accuracy, and consensus error) at 150 epochs, while \textit{Prox-DBRO-LSVRG} ranks second on these three performance metrics.
\begin{figure*}[!htp]
\begin{center}
\subfloat[A 60-agent network containing $20$ Byzantine agents.]{\includegraphics[width=1.7in,height=1.2in]{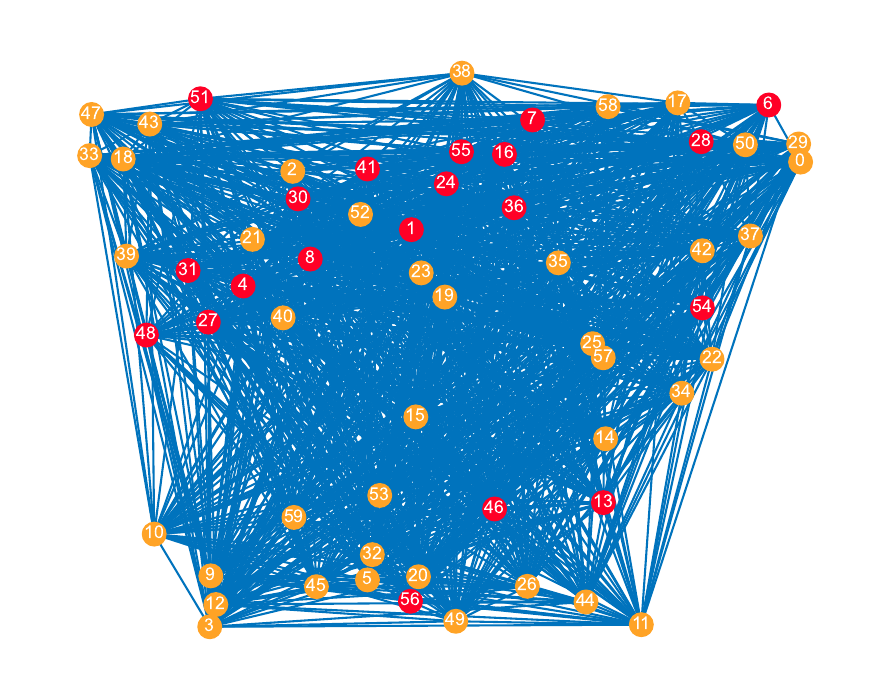}\label{Fig-3-1}} \hfill
\subfloat[Optimal gap over epochs.]{\includegraphics[width=1.7in,height=1.2in]{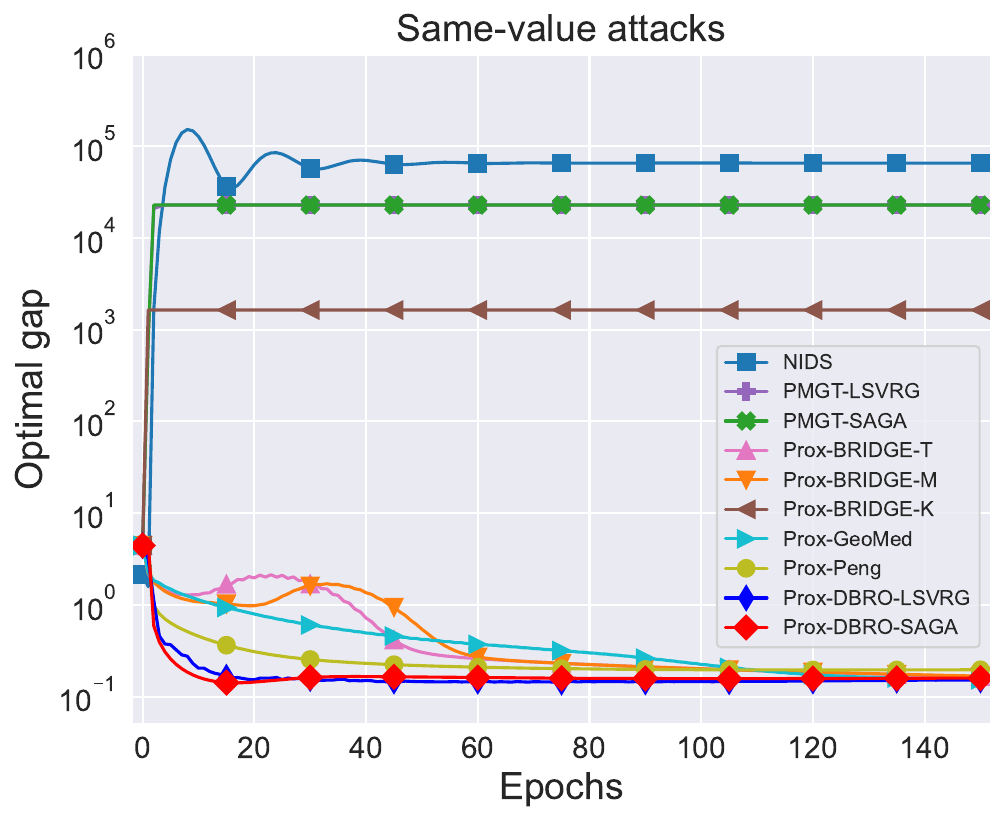}\label{Fig-3-2}} \hfill
\subfloat[Testing accuracy over epochs.]{\includegraphics[width=1.7in,height=1.2in]{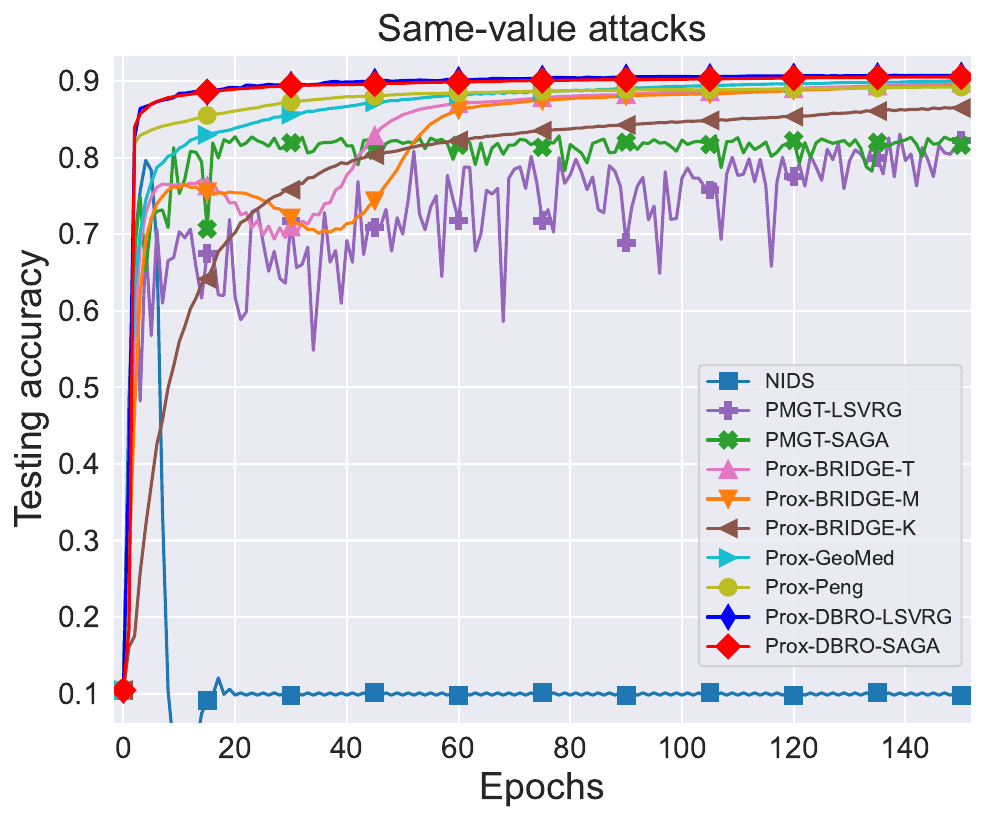}\label{Fig-3-3}} \hfill
\subfloat[Consensus error over epochs.]{\includegraphics[width=1.7in,height=1.2in]{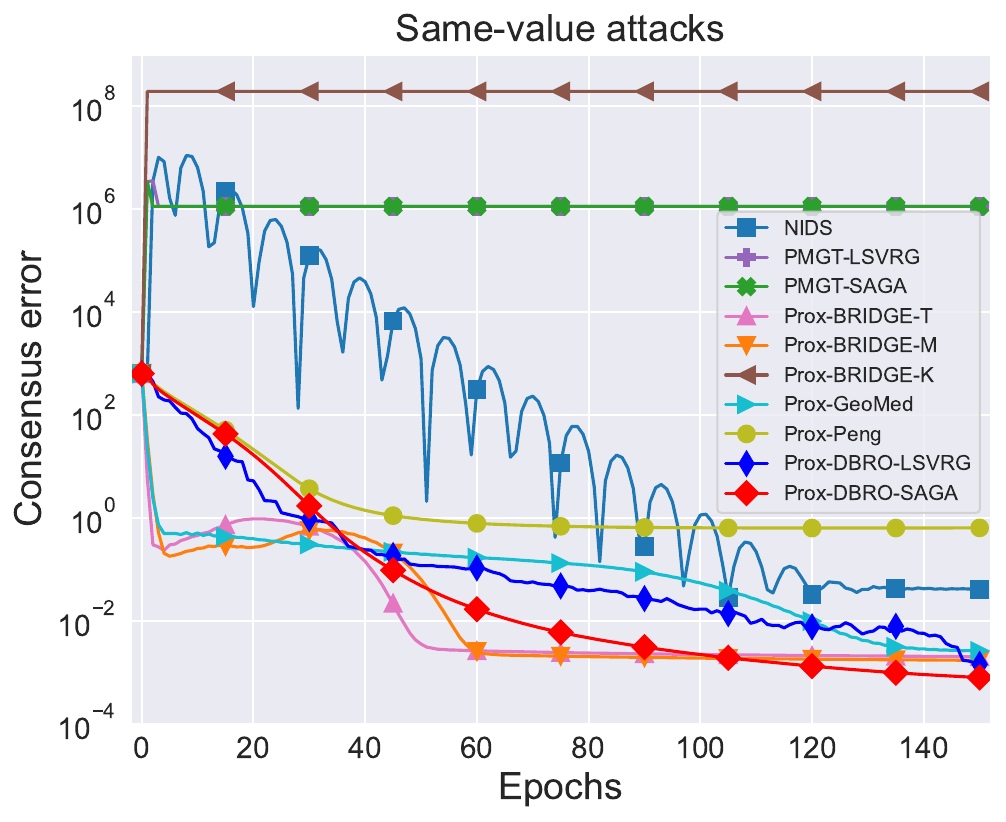}\label{Fig-3-4}} \hfill
\end{center}
\caption{Performance comparison among all tested algorithms under same-value attacks.}
\label{Fig-3}
\end{figure*}
\begin{table*}[!htp]
\centering
\caption{Parameter settings and algorithm performance at 150 epochs under same-value attacks.}
\scalebox{0.70}{
\begin{tabular}{ccccccccccc}
\toprule
\diagbox{PS\&PI}{Algorithms} & \textit{NIDS} & \textit{PMGT-LSVRG} & \textit{PMGT-SAGA} & \textit{Prox-BRIDGE-T} & \textit{Prox-BRIDGE-M} & \textit{Prox-BRIDGE-K} & \textit{Prox-GeoMed} & \textit{Prox-Peng} & {\bf{\textit{Prox-DBRO-LSVRG}}} & {\bf{\textit{Prox-DBRO-SAGA}}}\\
\midrule
Step-size & [0.5, 0.55] & 0.5 & 0.5 & 0.35 & 0.3 & 0.4 & 0.4 & $0.97/ \left({k + 25}\right)$ & $0.91/ \left({k + 21}\right) $ & $0.74/ \left({k + 35}\right) $  \\
Triggering probability & N/A & $m/Q$ & N/A & N/A & N/A & N/A & N/A & N/A & $\left[ {m/Q/8,m/Q/4} \right]$ & N/A  \\
Penalty parameter & N/A & N/A & N/A & N/A & N/A & N/A & N/A & 0.0005 & [0.0004, 0.00045] & [0.0005, 0.00055] \\
Consensus error & 4.1710e-02 & 1.1216e+06 & 1.1216e+06 & 2.0447e-03 & 1.7397e-03 & 1.9091e+08 & 2.6619e-03 & 6.0639e-01 & 1.4196e-03 & \bf{7.0507e-04} \\
Testing accuracy & 0.098 & 0.8227 & 0.8162 & 0.8972 & 0.8942 & 0.8649 & 0.9006 & 0.8915 & \bf{0.9072} & 0.9067 \\
Optimal gap & 6.5336e+04 & 2.2783e+04 & 2.2783e+04 & 1.6548e-01 & 1.6724e-01 & 1.6336e+03 & 1.5303e-01 & 1.9619e-01 & \bf{1.4623e-01} & 1.5938e-01 \\
\hline
\multicolumn{11}{c}{\makecell[c]{PS\&PI is the abbreviation of parameter settings and performance metrics.}}\\
\bottomrule
\end{tabular}}
\label{Tab-3}       
\end{table*}\newline
{\bf{Same-value attacks:}} As depicted in Fig. \ref{Fig-3}: \subref{Fig-3-1}, an $m=60$ network consists of ${\left| \mathcal{R} \right|}= 40$ reliable agents (yellow nodes) and $\left| \mathcal{B} \right| = 20$  Byzantine agents (red nodes), where each Byzantine agent $b$, $b \in \mathcal{B}_i$, keeps sending ${z_{ib,k}} = 1000 * {1_n}$ to its reliable neighbor $i$, $i \in {\mathcal{R}}$, at each iteration. Under this attack, the states of reliable agents can be easily blown up to sufficiently large values, which prevents the tested algorithms from convergence. Figs. \ref{Fig-3}: \subref{Fig-3-2}-\subref{Fig-3-3} manifest that the proposed algorithms achieve faster convergence than the other tested algorithms on the performance metrics of the optimal gap and testing accuracy, while \textit{Prox-DBRO-LSVRG} is slightly faster than \textit{Prox-DBRO-SAGA} in this case. It can be found in Table \ref{Tab-3} that the proposed algorithms attain also a smaller consensus error than the other tested algorithms at 150 epochs. Note that the consensus error of \textit{Prox-BRIDGE-B} goes to a very large value since it adopts a vector-wise screening at each reliable agent, which results in a single surviving vector totally from one neighboring agent and thus easily leads to a large state variance between any two reliable agents. The performance comparison takes no account of \textit{BRIDGE-B} \cite{Fang2022} due to its strict requirement on the number of neighboring agents and high computational overhead. In a nutshell, the above numerical experiments demonstrate that while the proposed algorithms do not achieve the smallest consensus error only under zero-sum attacks, they achieve the best performance in all other cases.
\begin{figure}[!htp]
\begin{center}
\subfloat[Optimal gap.]{\includegraphics[width=1.15in,height=0.8in]{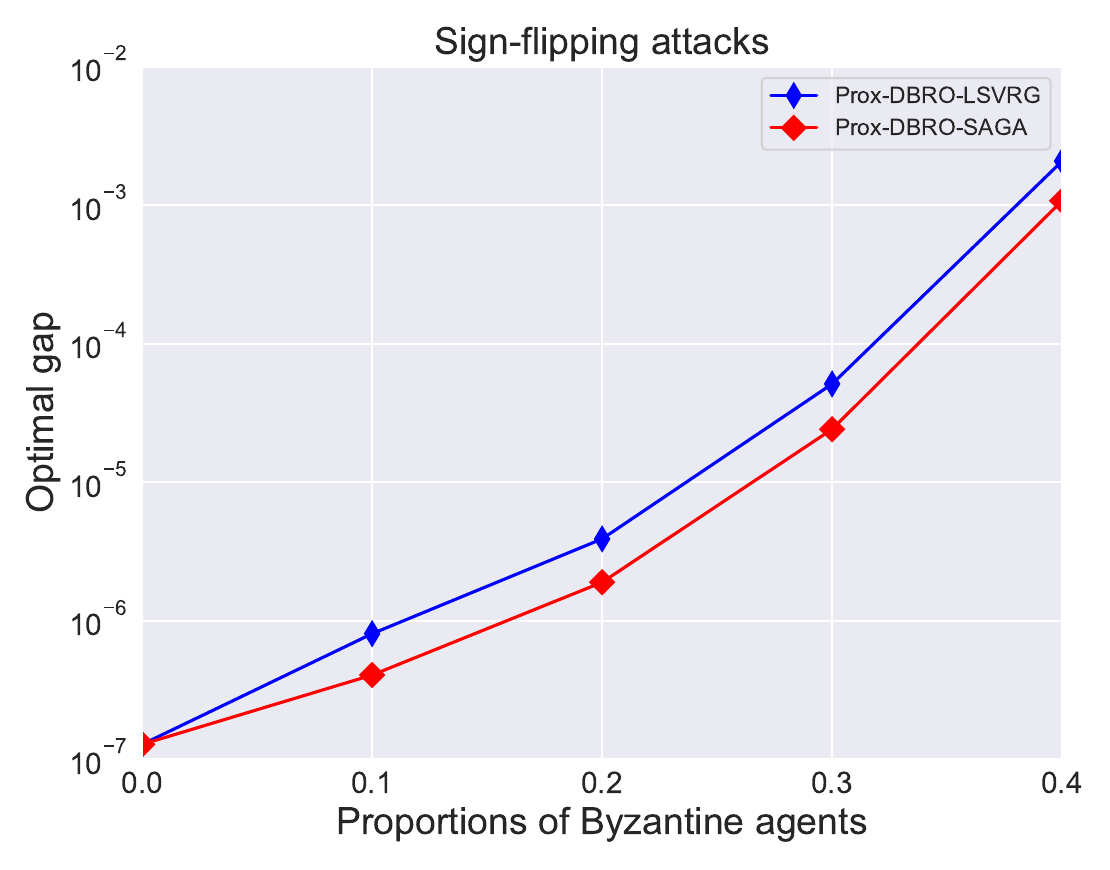}\label{Fig-4-1}} \hfill
\subfloat[Testing accuracy.]{\includegraphics[width=1.15in,height=0.8in]{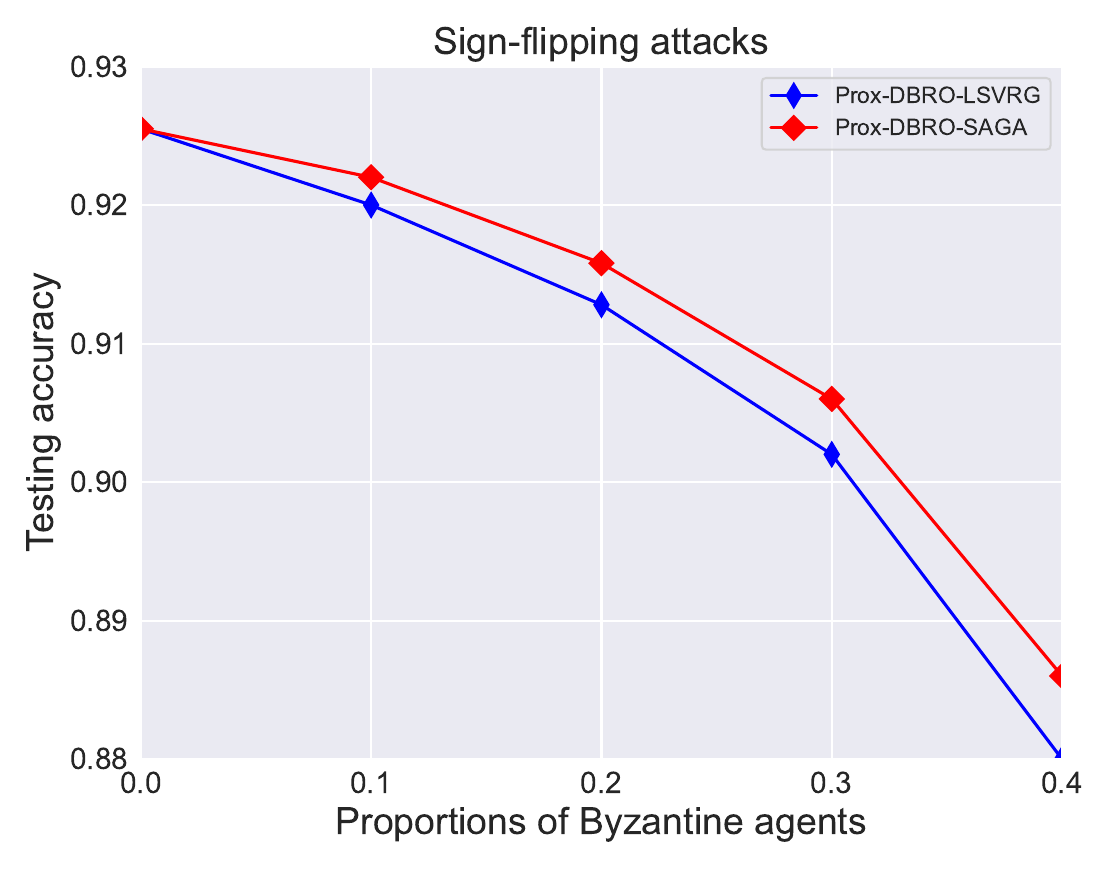}\label{Fig-4-2}} \hfill
\subfloat[Consensus error.]{\includegraphics[width=1.15in,height=0.8in]{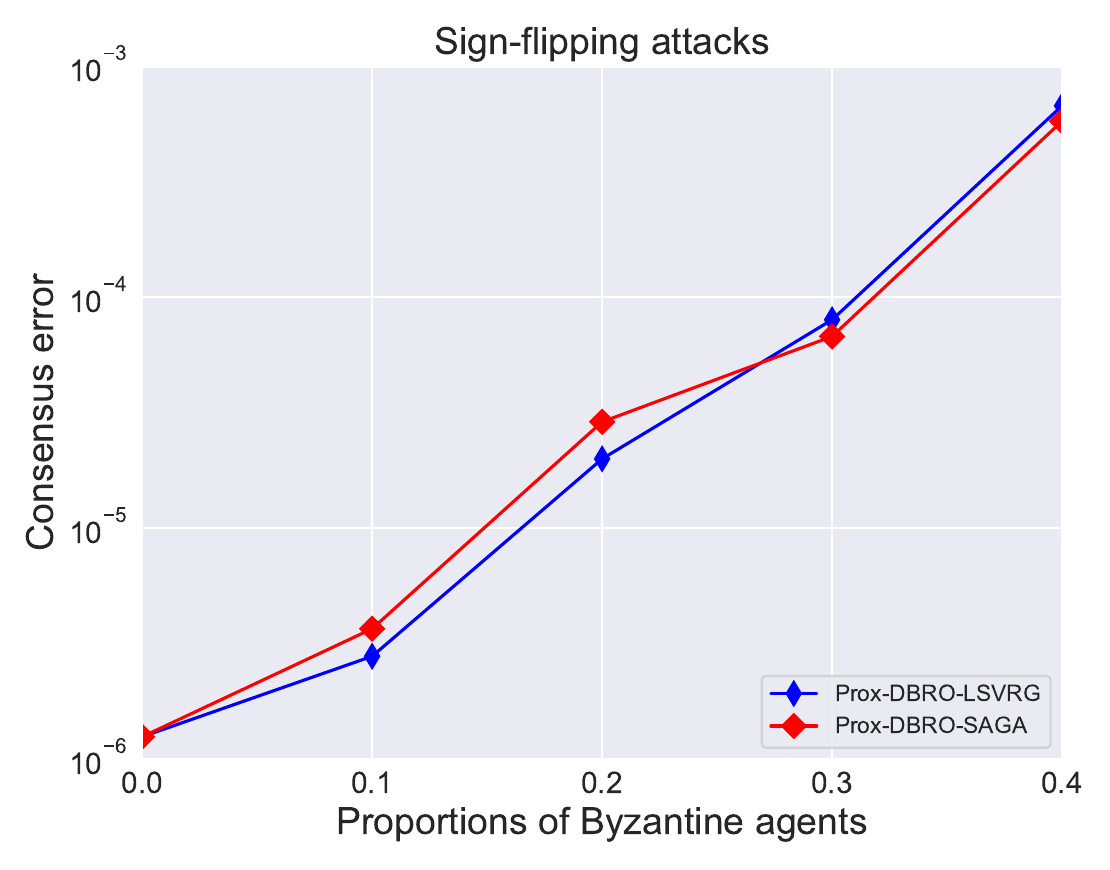}\label{Fig-4-3}} \hfill
\end{center}
\caption{Performance degradation versus increasing proportions of Byzantine agents under sign-flipping attacks.}
\label{Fig-4}
\end{figure}\newline
{\bf{Sign-flipping attacks:}} In this case, we aim to verify the trade-off between the convergence accuracy and Byzantine resilience of the propose algorithms.
The total number of agents including both reliable and Byzantine agents is $m = 100$, where each Byzantine agent $b$,  $b \in \mathcal{B}$, sends the falsified model ${z_{ib,k}} =  - s_b \sum\nolimits_{j \in {\mathcal{R}_i} \cup \left\{ i \right\}} {{x_{j,k}}} /\left( {\left| {{\mathcal{R}_i}} \right| + 1} \right)$ to their reliable neighbors $i$, $i \in {\mathcal{R}_b}$, where $s_b > 0$ is the hyperparameter controlling the deviation of the attack.
From Fig. \ref{Fig-4}, we can see that if the proportion or number of Byzantine agents increases, then the convergence accuracy regarding all three performance metrics becomes worse. This verifies the trade-off established in the theoretical result (see Theorems \ref{T2}-\ref{T3}).

\section{Conclusion}\label{Conclu}
In this paper, we proposed two Byzantine-resilient decentralized and VR stochastic gradient algorithms, namely \textit{Prox-DBRO-LSVRG} and \textit{Prox-DBRO-SAGA}, to resolve a category of nonsmooth composite finite-sum optimization problems over MASs in the presence of Byzantine agents. Theoretical analysis established both linear and sub-linear convergence rates for the proposed algorithms under different conditions. In the numerical experiments, the proposed algorithms were applied to resolving a decentralized sparse soft-max regression task over MASs under different Byzantine attacks, which verifies the theoretical findings and demonstrates the better convergence performance of the proposed algorithms than the other notable decentralized algorithms. However, both \textit{Prox-DBRO-LSVRG} and \textit{Prox-DBRO-SAGA} are not perfect, since they can only achieve exact sub-linear convergence in the absence of Byzantine agents. Future work will further investigate privacy issues and intermittent communication, which are also prevalent in MASs.

\appendix \label{Appendix}
\subsection{Proof of Lemma \ref{L1}}\label{Appen1}
According to Step 7 in Algorithm \ref{Algo3}, at iteration $k$, $\forall k \ge 1$, the auxiliary variables $u_{i,k + 1}^l$, $i \in \mathcal{R}$, take value $u_{i,k}^l$ or ${x_{i,k}}$, associated with probabilities $\left( {1 - 1/{q_i}} \right)$ and $1/q_i$, respectively.
This observation is owing to the fact that selection of the random sample for \textit{Prox-DBRO-SAGA}, at each iteration $k \ge 1$, is uniformly and independently executed. Hence, we have
\begin{equation}\label{E7-1}
\begin{aligned}
&{\mathbb{E}_k}\left[ {\frac{1}{{{q_i}}}\sum\limits_{l = 1}^{{q_i}} {{\nabla f_i^l{{\left( {{{\tilde x}^*}} \right)}^ \top }\left( {u_{i,k + 1}^l - {{\tilde x}^*}} \right)}} } \right]\\
=& \left( {1 - \frac{1}{{{q_i}}}} \right)\frac{1}{{{q_i}}}\sum\limits_{l = 1}^{{q_i}} {\nabla f_i^l{{\left( {{{\tilde x}^*}} \right)}^ \top }\left( {u_{i,k}^l - {{\tilde x}^*}} \right)} + \frac{1}{{{q_i}}}\nabla {f_i}{\left( {{{\tilde x}^*}} \right)^ \top }\\
& \times \left( {{x_{i,k}} - {{\tilde x}^*}} \right).
\end{aligned}
\end{equation}
Similarly, it holds that
\begin{equation}\label{E7-2}
{\mathbb{E}_k}\left[ {f_i^l\left( {u_{i,k + 1}^l} \right)} \right] = \left( {1 - \frac{1}{{{q_i}}}} \right)f_i^l\left( {u_{i,k}^l} \right) + \frac{1}{{{q_i}}}f_i^l\left( {{x_{i,k}}} \right).
\end{equation}
Via summing (\ref{E7-2}) over index $l$ for all $l = 1, \ldots ,{q_i}$, we can further obtain
\begin{equation}\label{E7-3}
\begin{aligned}
&{\mathbb{E}_k}\left[ {\frac{1}{{{q_i}}}\sum\limits_{l = 1}^{{q_i}} {f_i^l\left( {u_{i,k + 1}^l} \right)} } \right]\\
= & \left( {1 - \frac{1}{{{q_i}}}} \right)\frac{1}{{{q_i}}}\sum\limits_{l = 1}^{{q_i}} {f_i^l\left( {u_{i,k}^l} \right)}  + \frac{1}{{{q_i}}}\frac{1}{{{q_i}}}\sum\limits_{l = 1}^{{q_i}} {f_i^l\left( {{x_{i,k}}} \right)}\\
= & \left( {1 - \frac{1}{{{q_i}}}} \right)\frac{1}{{{q_i}}}\sum\limits_{l = 1}^{{q_i}} {f_i^l\left( {u_{i,k}^l} \right)}  + \frac{1}{{{q_i}}}{f_i}\left( {{x_{i,k}}} \right).
\end{aligned}
\end{equation}
Recalling the definition of $t_{i,k}^u$ and combining Eqs. (\ref{E7-1}) and (\ref{E7-3}), we obtain
\begin{equation}\label{E7-4}
\begin{aligned}
& {\mathbb{E}_k}\left[ {t_{i,k + 1}^u} \right]\\
=&{\mathbb{E}_k}\left[ {\frac{1}{{{q_i}}}\sum\limits_{l = 1}^{{q_i}} {f_i^l\left( {u_{i,k + 1}^l} \right)}  \!-\! {f_i}\left( {{{\tilde x}^*}} \right) \!-\! \nabla f_i^l{{\left( {{{\tilde x}^*}} \right)}^ \top }\left( {u_{i,k + 1}^l \!- {{\tilde x}^*}} \right)} \right]\\
=&\left( {1 - \frac{1}{{{q_i}}}} \right)\frac{1}{{{q_i}}}\sum\limits_{l = 1}^{{q_i}} {f_i^l\left( {u_{i,k}^l} \right) - \nabla f_i^l{{\left( {{{\tilde x}^*}} \right)}^ \top }\left( {u_{i,k}^l - {{\tilde x}^*}} \right)}\\
&+ \frac{1}{{{q_i}}}{f_i}\left( {{x_{i,k}}} \right) - {f_i}\left( {{{\tilde x}^*}} \right) - \frac{1}{{{q_i}}}\nabla {f_i}{\left( {{{\tilde x}^*}} \right)^ \top }\left( {{x_{i,k}} - {{\tilde x}^*}} \right)\\
=& \left( {1 - \frac{1}{{{q_i}}}} \right){t_{i,k}}  + \frac{1}{{{q_i}}}\left( {{f_i}\left( {{x_{i,k}}} \right) - {f_i}\left( {{{\tilde x}^*}} \right)} \right)- \frac{1}{{{q_i}}}\nabla {f_i}{\left( {{{\tilde x}^*}} \right)^ \top }\\
& \times\left( {{x_{i,k}} - {{\tilde x}^*}} \right).
\end{aligned}
\end{equation}
Summing Eq. (\ref{E7-4}) over $i$ yields
\begin{equation}\label{E7-5}
\begin{aligned}
&\sum\limits_{i \in \mathcal{R}} {{\mathbb{E}_k}\left[ {t_{i,k + 1}^u} \right]} \\
= &  \sum\limits_{i \in \mathcal{R}} {\frac{1}{{{q_i}}}\left( {{f_i}\left( {{x_{i,k}}} \right) - {f_i}\left( {{{\tilde x}^*}} \right) - \nabla {f_i}{{\left( {{{\tilde x}^*}} \right)}^ \top }\left( {{x_{i,k}} - {{\tilde x}^*}} \right)} \right)}\\
&+\sum\limits_{i \in \mathcal{R}} {\left( {1 - \frac{1}{{{q_i}}}} \right)t_{i,k}^u}\\
\le &\frac{1}{{{q_{\min }} }}\sum\limits_{i \in \mathcal{R}} { {{f_i}\left( {{x_{i,k}}} \right) - {f_i}\left( {{{\tilde x}^*}} \right) - \nabla {f_i}{{\left( {{{\tilde x}^*}} \right)}^ \top }\left( {{x_{i,k}} - {{\tilde x}^*}} \right)} }\\
&+\left( {1 - \frac{1}{{{q_{\max }}}}} \right)\sum\limits_{i \in \mathcal{R}} {t_{i,k}^u} \\
= & \frac{1}{{{q_{\min }} }}{D_F}\left( {{x_k},{x^*}} \right)+\left( {1 - \frac{1}{{{q_{\max }}}}} \right)\sum\limits_{i \in \mathcal{R}} {t_{i,k}^u},
\end{aligned}
\end{equation}
where the second inequality uses $1 \le {q_{\min }}  \le {q_i} \le {q_{\max }}$, and the last equality is in accordance with the definition of ${D_F}\left( {{x_k},{x^*}} \right)$. Substituting the definition of $t_k^u$ obtains the relation (\ref{E4-1}). In view of Step 7 in Algorithm \ref{Algo2}, we know that at iteration $k$, $\forall k \ge 1$, the auxiliary variables ${w_{i,k + 1}}$, $\forall i \in \mathcal{R}$, either take value ${x_{i,k}}$ with probability $p_i$, or keep the most recent update ${w_{i,k}}$ with probability $1 - {p_i}$. Therefore, it can be verified that
\begin{equation}\label{E7-6}
\begin{aligned}
&{\mathbb{E}_k}\left[ {\frac{1}{{{q_i}}}\sum\limits_{l = 1}^{{q_i}} { {\nabla f_i^l{{\left( {{{\tilde x}^*}} \right)}^ \top }\left( {{w_{i,k + 1}} - {{\tilde x}^*}} \right)} } } \right]\\
=& \frac{{{p_i}}}{{{q_i}}}\sum\limits_{l = 1}^{{q_i}} {\nabla f_i^l{{\left( {{{\tilde x}^*}} \right)}^ \top }\left( {{x_{i,k}} \!-\! {{\tilde x}^*}} \right)} \!+\! \left( {1 \!-\! {p_i}} \right)\nabla {f_i}{\left( {{{\tilde x}^*}} \right)^ \top }( {{w_{i,k}} \!-\! {{\tilde x}^*}} ).
\end{aligned}
\end{equation}
Likewise, we have
\begin{equation}\label{E7-7}
{\mathbb{E}_k}\!\!\left[ {\frac{1}{{{q_i}}}\sum\limits_{l = 1}^{{q_i}} {f_i^l\left( {{w_{i,k + 1}}} \right)} } \right]\\
\!\!\!=\! \left( {1 \!-\! {p_i}} \right)\frac{1}{{{q_i}}}\sum\limits_{l = 1}^{{q_i}} {f_i^l\left( {{w_{i,k}}} \right)}  + {p_i}{f_i}\left( {{x_{i,k}}} \right).
\end{equation}
Recalling the definition of $t_{i,k}^w$ and combining Eq. (\ref{E7-6}) with (\ref{E7-7}), we have
\begin{equation*}
\begin{aligned}
&{\mathbb{E}_k}\left[ {t_{i,k + 1}^w} \right]\\
=&{\mathbb{E}_k}\left[ {\frac{1}{{{q_i}}}\sum\limits_{l = 1}^{{q_i}} {f_i^l\left( {{w_{i,k + 1}}} \right)} \!-\! {f_i}\left( {{{\tilde x}^*}} \right) \! -\! \nabla f_i^l{{\left( {{{\tilde x}^*}} \right)}^ \top }\left( {{w_{i,k + 1}} \! -\!  {{\tilde x}^*}} \right)} \right]\\
\end{aligned}
\end{equation*}
\begin{equation}\label{E7-8}
\begin{aligned}
=& \left( {1 \!- {p_i}} \right){t_{i,k}^w} \!+\! {p_i}( {{f_i}\left( {{x_{i,k}}} \right) \!-\! {f_i}\left( {{{\tilde x}^*}} \right) \!-\! \nabla {f_i}{{\left( {{{\tilde x}^*}} \right)}^ \top }\left( {{x_{i,k}} \!-\! {{\tilde x}^*}} \right)} ),
\end{aligned}
\end{equation}
where we apply $ {f_i}\left( {{{\tilde x}^*}} \right) = \left( {1/{q_i}} \right)\sum\nolimits_{l = 1}^{{q_i}} {f_i^l\left( {{\tilde x}^*} \right)}$ in the last equality. The relation (\ref{E4-2}) is reached through summing Eq. (\ref{E7-8}) over $i$ and substituting the definitions of $t_k^w$ and ${D_F}\left( {{x_k},{x^*}} \right)$.

\subsection{Proof of Lemma \ref{L2}}\label{Appen2}
According to Step 3 in Algorithm \ref{Algo2}, it holds that
\begin{equation}\label{E8-1}
\begin{aligned}
&{\mathbb{E}_k}\left[ {\left\| {r_{i,k}^u - \nabla {f_i}\left( {{{\tilde x}^*}} \right)} \right\|_2^2} \right] \\
= & {\mathbb{E}_k}\left[ {\left\| {r_{i,k}^u - \nabla {f_i}\left( {{{\tilde x}^*}} \right) - \nabla {f_i}\left( {{x_{i,k}}} \right) + \nabla {f_i}\left( {{{\tilde x}^*}} \right)} \right\|_2^2} \right]\\
&+\left\| {\nabla {f_i}\left( {{x_{i,k}}} \right) - \nabla {f_i}\left( {{{\tilde x}^*}} \right)} \right\|_2^2,
\end{aligned}
\end{equation}
where the equality is due to the standard variance decomposition $\mathbb{E}_k\left[ {\left\| A \right\|_2^2} \right] = \left\| {\mathbb{E}_k\left[ A \right]} \right\|_2^2 + \mathbb{E}_k\left[ {\left\| {A - \mathbb{E}_k\left[ A \right]} \right\|_2^2} \right]$, with $A = {r_{i,k}^u} - \nabla {f_i}\left( {{{\tilde x}^*}} \right)$. We continue to handle the first term in the RHS of Eq. (\ref{E8-1}) as follows:
\begin{equation}\label{E8-2}
\begin{aligned}
&{\mathbb{E}_k}\left[ {\left\| {r_{i,k}^u - \nabla {f_i}\left( {{{\tilde x}^*}} \right) - \nabla {f_i}\left( {{x_{i,k}}} \right) + \nabla {f_i}\left( {{{\tilde x}^*}} \right)} \right\|_2^2} \right]\\
\le & 2{\mathbb{E}_k}\!\!\left[ {\left\| {\nabla f_i^{{s_{i,k}}}\!\left( {{x_{i,k}}} \right) \!-\! \nabla f_i^{{s_{i,k}}}\left( {{{\tilde x}^*}} \right) \!-\! \nabla {f_i}\left( {{x_{i,k}}} \right) \!+\! \nabla {f_i}\left( {{{\tilde x}^*}} \right)} \right\|_2^2} \right]\\
& + 2{\mathbb{E}_k}\!\left[ {\left\| {\nabla f_i^{{s_{i,k}}}\left( {u_{i,k}^{{s_{i,k}}}} \right) - \nabla f_i^{{s_{i,k}}}\left( {{{\tilde x}^*}} \right) - \frac{1}{{{q_i}}}\sum\limits_{l = 1}^{{q_i}} {\nabla f_i^l\left( {u_{i,k}^l} \right)} } \right.} \right.\\
&+ \left. {\left. {  \nabla {f_i}\left( {{{\tilde x}^*}} \right)} \right\|_2^2} \right]\\
\le & 2{\mathbb{E}_k}\!\!\left[ {\left\| {\nabla f_i^{{s_{i,k}}}\left( {{x_{i,k}}} \right) \!-\! \nabla f_i^{{s_{i,k}}}\left( {{{\tilde x}^*}} \right) \!-\! \nabla {f_i}\left( {{x_{i,k}}} \right) \!+\! \nabla {f_i}\left( {{{\tilde x}^*}} \right)} \right\|_2^2} \right]\\
&+ 2{\mathbb{E}_k}\left[ {\left\| {\nabla f_i^{{s_{i,k}}}\left( {u_{i,k}^{{s_{i,k}}}} \right) - \nabla f_i^{{s_{i,k}}}\left( {{{\tilde x}^*}} \right)} \right\|_2^2} \right]\\
=& 2{\mathbb{E}_k}\left[ {\left\| {\nabla f_i^{{s_{i,k}}}\left( {{x_{i,k}}} \right) - \nabla f_i^{{s_{i,k}}}\left( {{{\tilde x}^*}} \right)} \right\|_2^2} \right]\\
&+ 2{\mathbb{E}_k}\left[ {\left\| {\nabla f_i^{{s_{i,k}}}\left( {u_{i,k}^{{s_{i,k}}}} \right) - \nabla f_i^{{s_{i,k}}}\left( {{{\tilde x}^*}} \right)} \right\|_2^2} \right]\\
& - 2\left\| {\nabla {f_i}\left( {{x_{i,k}}} \right) - \nabla {f_i}\left( {{{\tilde x}^*}} \right)} \right\|_2^2,
\end{aligned}
\end{equation}
where the second inequality utilizes $\mathbb{E}_k\left[ {\left\| {B - \mathbb{E}_k\left[ B \right]} \right\|_2^2} \right] \le \mathbb{E}_k\left[ {\left\| B \right\|_2^2} \right]$, with $B = \nabla f_i^{{s_{i,k}}}\left( {u_{i,k}^{{s_{i,k}}}} \right) - \nabla f_i^{{s_{i,k}}}\left( {{{\tilde x}^*}} \right)$, and the last equality applies the standard variance decomposition again. We proceed to substitute (\ref{E8-2}) back into (\ref{E8-1}) to obtain
\begin{equation}\label{E8-3}
\begin{aligned}
&{\mathbb{E}_k}\left[ {\left\| {r_{i,k}^u - \nabla {f_i}\left( {{{\tilde x}^*}} \right)} \right\|_2^2} \right]\\
= & 2{\mathbb{E}_k}\left[ {\left\| {\nabla f_i^{{s_{i,k}}}\left( {{x_{i,k}}} \right) - \nabla f_i^{{s_{i,k}}}\left( {{{\tilde x}^*}} \right)} \right\|_2^2} \right]\\
&+ 2{\mathbb{E}_k}\left[ {\left\| {\nabla f_i^{{s_{i,k}}}\left( {u_{i,k}^{{s_{i,k}}}} \right) - \nabla f_i^{{s_{i,k}}}\left( {{{\tilde x}^*}} \right)} \right\|_2^2} \right]\\
&- \left\| {\nabla {f_i}\left( {{x_{i,k}}} \right) - \nabla {f_i}\left( {{{\tilde x}^*}} \right)} \right\|_2^2.
\end{aligned}
\end{equation}
Summing Eq. (\ref{E8-3}) over $i$ generates
\begin{equation}\label{E8-4}
\begin{aligned}
&{\mathbb{E}_k}\left[ {\left\| {r_k^u - \nabla F\left( {{{x}^*}} \right)} \right\|_2^2} \right]\\
\le & 2\sum\limits_{i \in \mathcal{R}} {{\mathbb{E}_k}\left[ {\left\| {\nabla f_i^{{s_{i,k}}}\left( {{x_{i,k}}} \right) - \nabla f_i^{{s_{i,k}}}\left( {{{\tilde x}^*}} \right)} \right\|_2^2} \right]} \\
&+ 2\sum\limits_{i \in \mathcal{R}} {{\mathbb{E}_k}\left[ {\left\| {\nabla f_i^{{s_{i,k}}}\left( {u_{i,k}^{{s_{i,k}}}} \right) - \nabla f_i^{{s_{i,k}}}\left( {{{\tilde x}^*}} \right)} \right\|_2^2} \right]} \\
&- \sum\limits_{i \in \mathcal{R}} {\left\| {\nabla {f_i}\left( {{x_{i,k}}} \right) - \nabla {f_i}\left( {{{\tilde x}^*}} \right)} \right\|_2^2}.
\end{aligned}
\end{equation}
Since the local component objective function $ f_i^l$, $\forall l \in \mathcal{Q}_i$ and $\forall i \in \mathcal{R}$, is $L$-smooth according to Assumption \ref{Assu1}, we have
\begin{equation}\label{E8-5}
\begin{aligned}
&\frac{1}{{2L}}\left\| {\nabla f_i^l\left( {u_{i,k}^l} \right) - \nabla {f_{i,l}}\left( {{{\tilde x}^*}} \right)} \right\|_2^2 \\
\le & f_i^l\left( {u_{i,k}^l} \right) - f_i^l\left( {{{\tilde x}^*}} \right) - \nabla f_i^l{\left( {{{\tilde x}^*}} \right)^ \top }\left( {u_{i,k}^l - {{\tilde x}^*}} \right).
\end{aligned}
\end{equation}
Summing the both sides of (\ref{E8-5}) over $l$ from $1$ to $q_i$ obtains
\begin{equation}\label{E8-6}
\frac{1}{{{q_i}}}\sum\limits_{l = 1}^{{q_i}} {\left\| {\nabla f_i^l\left( {u_{i,k}^l} \right) - \nabla f_i^l\left( {{{\tilde x}^*}} \right)} \right\|_2^2}  \le  2Lt_{i,k}^u.
\end{equation}
As the local component function $f_i^{{s_{i,k}}}$ has a uniform distribution over the set $\left\{ {f_i^1, \ldots ,f_i^{{q_i}}} \right\}$, it is natural to obtain
\begin{equation}\label{E8-7}
\begin{aligned}
&{\mathbb{E}_k}\left[ {\left\| {\nabla f_i^{{s_{i,k}}}\left( {u_{i,k}^{{s_{i,k}}}} \right) - \nabla f_i^{{s_{i,k}}}\left( {{{\tilde x}^*}} \right)} \right\|_2^2} \right]\\
 =& \frac{1}{{{q_i}}}\sum\limits_{l = 1}^{{q_i}} {\left\| {\nabla f_i^l\left( {u_{i,k}^l} \right) - \nabla f_i^l\left( {{{\tilde x}^*}} \right)} \right\|_2^2}.
 \end{aligned}
\end{equation}
Combining Eq. (\ref{E8-7}) and (\ref{E8-6}), and then summing over $i$, yields
\begin{equation}\label{E8-8}
\sum\limits_{i \in \mathcal{R}} {{\mathbb{E}_k}\left[ {\left\| {\nabla f_i^{{s_{i,k}}}\left( {u_{i,k}^{{s_{i,k}}}} \right) - \nabla f_i^{{s_{i,k}}}\left( {{{\tilde x}^*}} \right)} \right\|_2^2} \right]}  \le 2Lt_k^u.
\end{equation}
Summarizing (\ref{E8-4}) and (\ref{E8-8}) obtains
\begin{equation}\label{E8-9}
\begin{aligned}
&{\mathbb{E}_k}\left[ {\left\| {r_k^u - \nabla F\left( {{{x}^*}} \right)} \right\|_2^2} \right]\\
\le & 4Lt_k^u + 2\sum\limits_{i \in \mathcal{R}} {{\mathbb{E}_k}\left[ {\left\| {\nabla f_i^{{s_{i,k}}}\left( {{x_{i,k}}} \right) - \nabla f_i^{{s_{i,k}}}\left( {{{\tilde x}^*}} \right)} \right\|_2^2} \right]}\\
&- \left\| {\nabla F\left( {{x_k}} \right) - \nabla F\left( {{x^*}} \right)} \right\|_2^2,
\end{aligned}
\end{equation}
where we simplify $\sum\nolimits_{i \in \mathcal{R}} {\left\| {\nabla {f_i}\left( {{x_{i,k}}} \right) - \nabla {f_i}\left( {{{\tilde x}^*}} \right)} \right\|_2^2} $ as $\left\| {\nabla F\left( {{x_k}} \right) - \nabla F\left( {{x^*}} \right)} \right\|_2^2$. Via applying the Lipschitz continuity of $\nabla f_i^l$ again, we have
\begin{equation}\label{E8-10}
\sum\limits_{i \in \mathcal{R}} {{\mathbb{E}_k}\left[ {\left\| {\nabla f_i^{{s_{i,k}}}\left( {{x_{i,k}}} \right) - \nabla f_i^{{s_{i,k}}}\left( {{{\tilde x}^*}} \right)} \right\|_2^2} \right]} \le  2L{D_F}\left( {{x_k},{x^*}} \right),
\end{equation}
where we use the fact that $ {f_i}\left( {{x_{i,k}}} \right) = \left( {1/{q_i}} \right)\sum\nolimits_{l = 1}^{{q_i}} {f_i^l\left( {{x_{i,k}}} \right)}$ and $F\left( x_k \right) = \sum\nolimits_{i \in \mathcal{R}} {{f_i}\left( {{x_{i,k}}} \right)} $. Plugging (\ref{E8-10}) back into (\ref{E8-9}) generates
\begin{equation}\label{E8-11}
\begin{aligned}
&{\mathbb{E}_k}\left[ {\left\| {r_k^u - \nabla F\left( {{{x}^*}} \right)} \right\|_2^2} \right]\\
\le & 4Lt_k^u + 4L{D_F}\left( {{x_k},{x^*}} \right) - \left\| {\nabla F\left( {{x_k}} \right) - \nabla F\left( {{x^*}} \right)} \right\|_2^2.
\end{aligned}
\end{equation}
Considering $\mu$-strong convexity of the local objective function $f_i$, $\forall i \in \mathcal{R}$, we have
\begin{equation}\label{E8-12}
2\mu {D_F}\left( {{x_k},{x^*}} \right)\le \left\| {\nabla F\left( {{x_k}} \right) - \nabla F\left( {{x^*}} \right)} \right\|_2^2.
\end{equation}
Finally, one can obtain (\ref{E4-3}) via plugging the relation (\ref{E8-12}) back into (\ref{E8-11}). For \textit{Prox-DBRO-LSVRG}, we replace ${u_{i,k}^l}$ with ${{w_{i,k}}}$ to obtain (\ref{E4-4}), which completes the proof.

\subsection{Proof of Proposition \ref{P2}}\label{Appen3}
According to the definition of ${\mathbf{prox}}_{\alpha, G}\left\{ x \right\}$, we have
\begin{equation*}
\begin{aligned}
{\mathbf{prox}}_{\alpha, G}\left\{ x \right\} =& \arg \mathop {\min }\limits_y \left\{ {G\left( y \right) + \frac{1}{{2\alpha }}\left\| {y - x} \right\|_2^2} \right\}\\
=& \arg \mathop {\min }\limits_y \left\{ {\sum\limits_{i \in \mathcal{R}} {g\left( {{y_i}} \right)}  + \frac{1}{{2\alpha }}\sum\limits_{i \in \mathcal{R}} {\left\| {{y_i} - {x_i}} \right\|_2^2} } \right\}\\
\end{aligned}
\end{equation*}
\begin{equation}\label{E8-13}
\begin{aligned}
=&\left( {\begin{array}{*{20}{c}}
  {\arg \mathop {\min }\limits_{\tilde y \in {\mathbb{R}^n}} \left\{ {g\left( {\tilde y} \right) + \frac{1}{{2\alpha }}\left\| {\tilde y - {x_1}} \right\|_2^2} \right\}} \\
   \vdots  \\
  {\arg \mathop {\min }\limits_{\tilde y \in {\mathbb{R}^n}} \left\{ {g\left( {\tilde y} \right) + \frac{1}{{2\alpha }}\left\| {\tilde y - {x_{\left| \mathcal{R} \right|}}} \right\|_2^2} \right\}}
\end{array}} \right),
\end{aligned}
\end{equation}
which indicates ${\left[{\mathbf{prox}}_{\alpha, G}\left\{ x \right\}\right]}_i = {\mathbf{prox}}_{\alpha ,g}\left\{ x_i \right\}$. Based on this equality, it is straightforward to verify (\ref{E4-5}) with the aid of the non-expansiveness of the proximal operator ${{\mathbf{prox}}_{\alpha ,g}}$, which completes the proof.

\subsection{Proof of Theorem \ref{T1}}\label{Appen4}
The optimal solution to (\ref{E2-4}) satisfies the optimality condition
\begin{equation}\label{E9-2}
{0_n} \in {\nabla {f_i}\left( {x_i^*} \right) + \partial g\left( {x_i^*} \right)} + \frac{{{\phi}}}{2}\sum\limits_{j \in {\mathcal{R}_i}} {\partial {{\left\| {x_i^* - x_j^*} \right\|}_a}}, \forall i \in \mathcal{R}.
\end{equation}
According to the definition of the subdifferential ${\partial}{\left\| {x_i^* - x_j^*} \right\|_a} = \left\{ {{y_{ij}} \in {\mathbb{R}^n}|\left\langle {{y_{ij}},x_i^*} \right\rangle  = {{\left\| {x_i^*} \right\|}_a},{{\left\| {{y_{ij}}} \right\|}_b} \le 1} \right\}$, there exist $g'\left( {x_i^*} \right) \in \partial g\left( {x_i^*} \right)$ and ${\tilde y_{ij}} \in {\partial}{\left\| {x_i^* - x_j^*} \right\|_a}$, such that for $ i \in \mathcal{R}$
\begin{equation}\label{E9-3}
\begin{aligned}
{\nabla {f_i}\left( {x_i^*} \right) + g'\left( {x_i^*} \right)} + {\phi}\left( {\sum\limits_{j \in {\mathcal{R}_i},i < j} {{{\tilde y}_{ij}}}  - \sum\limits_{j \in {\mathcal{R}_i},i > j} {{{\tilde y}_{ji}}} } \right) = 0_n.
\end{aligned}
\end{equation}
Under Assumption \ref{Assu1}, the globally optimal solution $x^*$ exists uniquely. We next need to prove that the optimal solution ${\tilde x^*}$ satisfies (\ref{E9-3}), such that
\begin{equation}\label{E9-4}
{\nabla {f_i}\left( {{{\tilde x}^*}} \right) + g'\left( {{{\tilde x}^*}} \right)} + {\phi}\left( {\sum\limits_{j \in {\mathcal{R}_i},i < j} {{{\tilde y}_{ij}}}  - \sum\limits_{j \in {\mathcal{R}_i},i > j} {{{ \tilde y}_{ji}}} }  \right) = 0_n,
\end{equation}
where $g'\left( {{{\tilde x}^*}} \right) \in {\partial}g\left( {{{\tilde x}^*}} \right)$.
Since (\ref{E9-4}) admits an element-wise decomposition, hence without loss of generality, the rest proof assumes $n=1$, i.e., the scalar case. Via denoting $\Psi : = {\mathbf{col}}{\left\{ {{\psi _i}} \right\}_{i \in \mathcal{R}}} \in {\mathbb{R}^{\left| \mathcal{R} \right|}}$ with ${\psi _i}: = {\nabla {f_i}\left( {{{\tilde x}^*}} \right) + g'\left( {{{\tilde x}^*}} \right)} $, the task to prove (\ref{E9-4}) reduces to solving for a vector ${\tilde y}$, such that the following relation holds
\begin{subequations}\label{E9-5}
\begin{align}
\label{E9-5-1}&\phi \Pi  {\tilde y} + \Psi = 0_{\left| \mathcal{R} \right|},\\
\label{E9-5-2}&{\left\| {{{\tilde y}}} \right\|_b} \le  1,
\end{align}
\end{subequations}
where $\tilde y \in {\mathbb{R}^{\left| {{\mathcal{E}_\mathcal{R}}} \right|}}$ is the collected form of $\tilde y_{ij}$ according to the order of edges in ${\mathcal{E}_\mathcal{R}}$. We need to solve for at least one solution $\tilde y$ meeting ${\left\| {{{\tilde y}}} \right\|_b} \le 1$ with $b \ge 1$, such that (\ref{E9-5-1}) holds true. To proceed, we decompose the task into two parts. \\
{\bf{Part I:}} We first manifest that (\ref{E9-5-1}) has at least one solution. In view of the rank of the node-edge incidence matrix $\Pi$ is $\left| \mathcal{R} \right| - 1$ and the null space of the columns is spanned by the all-one vector ${1_{\left| \mathcal{R} \right|}}$. Recalling the definition of ${\psi _i}$, the optimality condition of (\ref{E2-1}) is $\sum\nolimits_{i \in \mathcal{R}} {{\psi _i}}  = 0$. Therefore, we know that the columns of $\Pi $ and those of $\left[ {\phi \Pi, \Psi } \right]$ share the same null space, which indicates the same rank of $\Pi$ and $\left[ {\phi \Pi , \Psi} \right]$ according to the rank-nullity theorem. The existence of solutions to (\ref{E9-5}) can be demonstrated according to the property of non-homogeneous linear equations. \\
{\bf{Part II:}} We provide a specific condition to ensure that the solution to (\ref{E9-5-1}) satisfies (\ref{E9-5-2}) as well. Suppose that $y \in {\mathbb{R}^{\left| {{\mathcal{E}_\mathcal{R}}} \right|}}$ is a solution to (\ref{E9-5-1}), such that $\phi \Pi  {y} + \Psi = 0_{\left| \mathcal{R} \right|}$.
 We consider the least-squares solution $y =  - {\Pi ^\dag }\Psi/\phi$, where $ {\Pi ^\dag }$ is the Moore-Penrose pseudo-inverse of $\Pi$. Then, it suffices to prove that ${\left\| y \right\|_b} \le 1$. Since ${\left\| y \right\|_b} = {\left( {\sum\nolimits_{i = 1}^{\left| \mathcal{R} \right|} {{{\left| {{y_i}} \right|}^b}} } \right)^{1/b}},\forall b \ge 1$, we know that ${\left\| y \right\|_b} \le {\left\| y \right\|_1}$. Therefore, we derive
\begin{equation}\label{E9-6}
\begin{aligned}
{\left\| y \right\|_b} \le & \frac{1}{\phi }{\left\| {{\Pi ^\dag }\Psi } \right\|_1}\\
\le & \frac{1}{\phi }{\left\| {{\Pi^\dag }} \right\|_1}{\left\| \Psi \right\|_1}\\
\le & \frac{1}{\phi } {\left| \mathcal{R} \right|\sqrt {\left| {{\mathcal{E}_\mathcal{R}}} \right|} }{\left\| {{\Pi^\dag }} \right\|_2}{\left\| \Psi \right\|_2},
\end{aligned}
\end{equation}
where the second inequality uses the vector-matrix norm compatibility, and the last inequality applies the facts that ${\left\| \Psi \right\|_1} \le \left| \mathcal{R} \right|{\left\| \Psi \right\|_2}$ and ${\left\| {{\Pi^\dag }} \right\|_1} \le \sqrt {\left| {{\mathcal{E}_\mathcal{R}}} \right|} {\left\| {{\Pi^\dag }} \right\|_2}$.
Consider ${{{\lambda }_{\max }}}\left( {{\Pi^\dag }} \right)$ and ${{{\lambda }_{\min }}\left( \Pi \right)}$ as the maximum and minimum singular values of matrices $\Pi^\dag$ and $\Pi$, respectively. Based on (\ref{E9-6}), we further obtain
\begin{equation}\label{E9-7}
\begin{aligned}
{\left\| y \right\|_b} \le & {{{\lambda }_{\max }}}\left( {{\Pi^\dag }} \right)\frac{{{{\left| \mathcal{R} \right|\sqrt {\left| {{\mathcal{E}_\mathcal{R}}} \right|} }}}}{{{\phi}}}{\left\| \Psi  \right\|_2}\\
=& \frac{{{{\left| \mathcal{R} \right|\sqrt {\left| {{\mathcal{E}_\mathcal{R}}} \right|} }}}}{{{\phi}{{{\lambda }_{\min }}\left( \Pi \right)}}}{\left\| \Psi \right\|_2}.
\end{aligned}
\end{equation}
Since ${\left\| \Psi \right\|_2} \le \sqrt {\left| \mathcal{R} \right|} {\left\| \Psi \right\|_\infty }$, we further have
\begin{equation}\label{E9-8}
{\left\| y \right\|_b} \le \frac{{{{\left| \mathcal{R} \right|}^{\frac{3}{2}}}\sqrt {\left| {{\mathcal{E}_\mathcal{R}}} \right|} }}{{{{{\lambda }_{\min }}\left( \Pi \right)}{\phi}}}\mathop {\max }\limits_{i \in \mathcal{R}} \left| {{\psi_i}} \right|.
\end{equation}
If we consider $n \ge 1$, i.e, the arbitrary dimension case, (\ref{E9-8}) becomes
\begin{equation}\label{E9-9}
{\left\| y \right\|_b} \le \frac{{{{\left| \mathcal{R} \right|}^{\frac{3}{2}}}\sqrt {\left| {{\mathcal{E}_\mathcal{R}}} \right|} }}{{{{{\lambda }_{\min }}\left( \Pi \right)}{\phi}}}\mathop {\max }\limits_{i \in \mathcal{R}} {\left\| \nabla {f_i}\left( {\tilde x^*} \right) + {g'\left( {{{\tilde x}^*}} \right)}\right\|_\infty }.
\end{equation}
The proof is completed by choosing an appropriate ${\phi}$ to meet
\begin{equation*}
\frac{{{{\left| \mathcal{R} \right|}^{\frac{3}{2}}}\sqrt {\left| {{\mathcal{E}_\mathcal{R}}} \right|} }}{{{{{\lambda }_{\min }}\left( \Pi \right)}{\phi}}}\mathop {\max }\limits_{i \in \mathcal{R}} {\left\| \nabla {f_i}\left( {\tilde x^*} \right) + {g'\left( {{{\tilde x}^*}} \right)} \right\|_\infty } \le 1,
\end{equation*}
which provides the condition to ensure the validity of (\ref{E9-5}) and thus finishes the proof.
\subsection{Proof of Theorem \ref{T2}}\label{Appen5}
Without loss of generality, we let ${{\tilde y}_{ji}} =  - {{\tilde y}_{ij}}$ for $j > i$ such that by defining ${{\tilde y}_i}: = \sum\nolimits_{j \in {{\mathcal R}_i}} {{{\tilde y}_{ij}}} $ and ${\partial _x}\chi \left( {{x^*}} \right): = {\mathbf{col}}{\left\{ {{{\tilde y}_i}} \right\}_{i \in \mathcal{R}}} \in {\mathbb{R}^{\left| \mathcal{R} \right|n}}$, it follows from the optimality condition (\ref{E9-2}) that
\begin{equation}\label{E10-1}
{x^*} = {\mathbf{prox}}_{\alpha, G}\left\{ {{x^*} - \alpha \left( {\nabla F\left( {{x^*}} \right) + {\partial _x}\chi \left( {{x^*}} \right)} \right)} \right\}.
\end{equation}
In view of the compact form (\ref{E3-5}), it holds
\begin{equation*}
\begin{aligned}
&{\mathbb{E}_k}\left[ {\left\| {{x_{k + 1}} - {x^*}} \right\|_2^2} \right]\\
=&{\mathbb{E}_k}\left[ {\left\| {{\mathbf{pro}}{{\mathbf{x}}_{\alpha ,G}}\left\{ {{{\bar x}_k}} \right\}-{\mathbf{pro}}{{\mathbf{x}}_{\alpha ,G}}\left\{ {{x^*}} \right.} \right. - \alpha \nabla F\left( {{x^*}} \right)} \right.\\
&- \left. {\left. {\left. {  \alpha {\partial _x}\chi \left( {{x^*}} \right)} \right\}} \right\|_2^2} \right]\\
\le & {\mathbb{E}_k}\left[ {\left\| {{{{\bar x}_k}} - \left( {{x^*} - \alpha \left( {\nabla F\left( {{x^*}} \right) + {\partial _x}\chi \left( {{x^*}} \right)} \right)} \right)} \right\|_2^2} \right]\\
=& \left\| {{x_k} - {x^*}} \right\|_2^2 - 2\alpha {\mathbb{E}_k}\left[ {\left\langle {{x_k} - {x^*},r_k - \nabla F\left( {{x^*}} \right)} \right\rangle } \right]\\
& - 2\alpha \left\langle {{x_k} - {x^*},{\partial _x}\chi \left( {{x_k}} \right) - {\partial _x}\chi \left( {{x^*}} \right) + {\partial _x}\delta \left( {{x_k}} \right)} \right\rangle \\
&+ {\alpha ^2}{\mathbb{E}_k}\left[ {\left\| {r_k - \nabla F\left( {{x^*}} \right) + {\partial _x}\chi \left( {{x_k}} \right) - {\partial _x}\chi \left( {{x^*}} \right)} \right.} \right.\\
\end{aligned}
\end{equation*}
\begin{equation}\label{E10-2}
\begin{aligned}
& + \left. {\left. {  {\partial _x}\delta \left( {{x_k}} \right)} \right\|_2^2} \right],
\end{aligned}
\end{equation}
where the second inequality applies (\ref{E4-5}). We continue to seek an upper bound for ${\mathbb{E}_k}\left[ {\left\| {r_k - \nabla F\left( {{x^*}} \right) + {\partial _x}\chi \left( {{x_k}} \right) - {\partial _x}\chi \left( {{x^*}} \right) + {\partial _x}\delta \left( {{x_k}} \right)} \right\|_2^2} \right]$ in the sequel.
\begin{equation}\label{E10-3}
\begin{aligned}
&{\mathbb{E}_k}\left[ {\left\| {r_k - \nabla F\left( {{x^*}} \right) + {\partial _x}\chi \left( {{x_k}} \right) - {\partial _x}\chi \left( {{x^*}} \right) + {\partial _x}\delta \left( {{x_k}} \right)} \right\|_2^2} \right]\\
\le & 4{\mathbb{E}_k}\left[ {\left\| {r_k - \nabla F\left( {{x^*}} \right)} \right\|_2^2} \right] + 2\left\| {{\partial _x}\chi \left( {{x_k}} \right) - {\partial _x}\chi \left( {{x^*}} \right)} \right\|_2^2\\
& + 4\left\| {{\partial _x}\delta \left( {{x_k}} \right)} \right\|_2^2\\
\le & 4\left( {4Lt_k + 2\left( {2L - \mu } \right){D_F}\left( {{x_k},{x^*}} \right)} \right) + 4\left\| {{\partial _x}\delta \left( {{x_k}} \right)} \right\|_2^2\\
& + 2\left\| {{\partial _x}\chi \left( {{x_k}} \right) - {\partial _x}\chi \left( {{x^*}} \right)} \right\|_2^2,\\
\end{aligned}
\end{equation}
where the first inequality applies ${\left\| {c + d} \right\|^2} \le 2{c^2} + 2{d^2}$ twice, and the second equality uses Lemma \ref{L2}. To proceed, we bound $\left\| {{\partial _x}\delta \left( {{x_k}} \right)} \right\|_2^2$ as follows:
\begin{equation}\label{E10-4}
\begin{aligned}
\left\| {{\partial _x}\delta \left( {{x_k}} \right)} \right\|_2^2=& \sum\limits_{i \in \mathcal{R}} {\left\| {{\phi}\sum\limits_{j \in {\mathcal{B}_i}} {{\partial _{{x_i}}}{{\left\| {{x_{i,k}} - {z_{ij,k}}} \right\|}_a}} } \right\|} _2^2\\
\le& n{{\phi }^2}\sum\limits_{i \in \mathcal{R}} {{{\left| {{\mathcal{B}_i}} \right|}^2}},
\end{aligned}
\end{equation}
where the inequality holds true since the $b$-norm ($b \ge 1$) of ${\partial _{{x_i}}}{\left\| {{x_{i,k}} - {z_{ij,k}}} \right\|_a}$, $\forall i \in \mathcal{R}$, is no larger than $1$ owing to Proposition \ref{P1}, i.e.,
\begin{equation}\label{E10-5}
\left| {{{\left[ {{\partial _{{x_i}}}{{\left\| {{x_{i,k}} - {z_{ij,k}}} \right\|}_a}} \right]}_e}} \right| \le 1,\forall e = 1, \ldots ,n.
\end{equation}
Following the same technical line as in (\ref{E10-4})-(\ref{E10-5}), it is not difficult to verify
\begin{equation}\label{E10-6}
\begin{aligned}
&\left\| {{\partial _x}\chi \left( {{x_k}} \right) - {\partial _x}\chi \left( {{x^*}} \right)} \right\|_2^2\\
=&\sum\limits_{i \in \mathcal{R}} {\left\| {{\phi}\sum\limits_{j \in {\mathcal{R}_i}} {\left( {{\partial _{{x_i}}}{{\left\| {{x_{i,k}} - {x_{j,k}}} \right\|}_a} - {\partial _{{x_i}}}{{\left\| {x_i^* - x_j^*} \right\|}_a}} \right)} } \right\|_2^2} \\
\le & 4n{{\phi }^2}\sum\limits_{i \in \mathcal{R}} {{{\left| {{\mathcal{R}_i}} \right|}^2}}.
\end{aligned}
\end{equation}
Combining (\ref{E10-3}), (\ref{E10-4}), and (\ref{E10-6}) obtains
\begin{equation}\label{E10-7}
\begin{aligned}
&{\mathbb{E}_k}\left[ {\left\| {r_k - \nabla F\left( {{x^*}} \right) + {\partial _x}\chi \left( {{x_k}} \right) - {\partial _x}\chi \left( {{x^*}} \right) + {\partial _x}\delta \left( {{x_k}} \right)} \right\|_2^2} \right]\\
\le & 4\left( {2Lt_k + \left( {2L - \mu } \right){D_F}\left( {{x_k},{x^*}} \right)} \right) + 8n{{\phi }^2}\sum\limits_{i \in \mathcal{R}} {{{\left| {{\mathcal{R}_i}} \right|}^2}}\\
& + 4n{{\phi }^2}\sum\limits_{i \in \mathcal{R}} {{{\left| {{\mathcal{B}_i}} \right|}^2}}.
\end{aligned}
\end{equation}
Since the local objective functions ${f_i}\left( {{x_i}} \right) $, $\forall i \in \mathcal{R}$, are $\mu$-strongly convex and $L$-smooth according to Assumption \ref{Assu1}, we have
\begin{equation}\label{E10-8}
\begin{aligned}
&- {\mathbb{E}_k}\left[ {\left\langle {{x_k} - {x^*},r_k - \nabla F\left( {{x^*}} \right)} \right\rangle } \right]\\
=& - \left\langle {{x_k} - {x^*},\nabla F\left( {{x_k}} \right) - \nabla F\left( {{x^*}} \right)} \right\rangle   \\
\le &- \frac{{\mu L}}{{\mu  + L}}\left\| {{x_k} - {x^*}} \right\|_2^2 - \frac{{1}}{{\mu  + L}}\left\| {\nabla F\left( {{x_k}} \right) - \nabla F\left( {{x^*}} \right)} \right\|_2^2.
\end{aligned}
\end{equation}
Recalling the definition of $\chi \left( {{x_k}} \right)$, we know that it is a convex function. Therefore, it is straightforward to obtain
\begin{equation}\label{E10-9}
- \left\langle {{x_k} - {x^*},{\partial _x}\chi \left( {{x_k}} \right) - {\partial _x}\chi \left( {{x^*}} \right)} \right\rangle  \le 0.
\end{equation}
Applying the Young's inequality and the relation (\ref{E10-4}) to the term $ - 2\left\langle {{x_k} - {x^*},\partial_x \delta \left( {{x_k}} \right)} \right\rangle $ yields
\begin{equation}\label{E10-10}
- 2\left\langle {{x_k} - {x^*},\partial_x \delta \left( {{x_k}} \right)} \right\rangle \le \gamma \left\| {{x_k} - {x^*}} \right\|_2^2 + \frac{{n{{\phi }^2}}}{\gamma }\sum\limits_{i \in \mathcal{R}} {{{\left| {{\mathcal{B}_i}} \right|}^2}},
\end{equation}
where the constant $\gamma  > 0$. Plugging the results (\ref{E10-7})-(\ref{E10-10}) into (\ref{E10-2}) gives
\begin{equation}\label{E10-11}
\begin{aligned}
&{\mathbb{E}_k}\left[ {\left\| {{x_{k + 1}} - {x^*}} \right\|_2^2} \right]\\
\le & \left( {1 - \left( {\frac{{2\mu L}}{{\mu  + L}} - \gamma } \right)\alpha } \right)\left\| {{x_k} - {x^*}} \right\|_2^2 +  8n{{\phi }^2}{\alpha ^2}\sum\limits_{i \in \mathcal{R}} {{{\left| {{\mathcal{R}_i}} \right|}^2}} \\
&+ 4{\alpha ^2}\left( {2Lt_k + \left( {2L - \mu } \right){D_F}\left( {{x_k},{x^*}} \right)} \right) + 4n{{\phi }^2}{\alpha ^2}\sum\limits_{i \in \mathcal{R}} {{{\left| {{\mathcal{B}_i}} \right|}^2}} \\
&+\frac{{n\alpha }}{\gamma }{{\phi }^2}\sum\limits_{i \in \mathcal{R}} {{{\left| {{\mathcal{B}_i}} \right|}^2}}.
\end{aligned}
\end{equation}
By setting $\gamma  = \mu L/\left( {\mu  + L} \right)$, we can rewrite (\ref{E10-11}) as follows:
\begin{equation}\label{E10-12}
\begin{aligned}
&{\mathbb{E}_k}\left[ {\left\| {{x_{k + 1}} - {x^*}} \right\|_2^2} \right] \\
\le & \left( {1 - \gamma {\alpha }} \right)\left\| {{x_k} \!-\! {x^*}} \right\|_2^2 \!+\! 4{\alpha ^2}\left( {2Lt_k + \left( {2L \!-\! \mu } \right){D_F}\left( {{x_k},{x^*}} \right)} \right)\\
&+ 8n{{\phi }^2}{\alpha ^2}\sum\limits_{i \in \mathcal{R}} {{{\left| {{\mathcal{R}_i}} \right|}^2}}+ 4n{{\phi }^2}{\alpha ^2}\sum\limits_{i \in \mathcal{R}} {{{\left| {{\mathcal{B}_i}} \right|}^2}}\!+\!\frac{{n\alpha }}{\gamma }{{\phi }^2}\sum\limits_{i \in \mathcal{R}} {{{\left| {{\mathcal{B}_i}} \right|}^2}}.
\end{aligned}
\end{equation}
According to (\ref{E4-1}), we have for any $c > 0$,
\begin{equation}\label{E10-13}
c\left( {{\mathbb{E}_k}\left[ {t_{k + 1}} \right] - t_k} \right) \le  - \frac{c}{{{q_{\max }}}}t_k + \frac{c}{{{q_{\min }}}}{D_F}\left( {{x_k},{x^*}} \right).
\end{equation}
Recall the definitions of $P_1$ and ${P_2}$. Combining (\ref{E10-12}) and (\ref{E10-13}) yields
\begin{equation}\label{E10-14}
\begin{aligned}
&{\mathbb{E}_k}\left[ {\left\| {{x_{k + 1}} - {x^*}} \right\|_2^2} \right] + c\left( {{\mathbb{E}_k}\left[ {t_{k + 1}} \right] - t_k} \right)\\
\le & \left( {1 - \gamma {\alpha}} \right)\left\| {{x_k} - {x^*}} \right\|_2^2 + 4n{{\phi }^2}\alpha^2\left(  {2\sum\limits_{i \in \mathcal{R}} {{{\left| {{\mathcal{R}_i}} \right|}^2}}  + \sum\limits_{i \in \mathcal{R}} {{{\left| {{\mathcal{B}_i}} \right|}^2}} }  \right)\\
&+ \frac{{n{{\phi }^2}}}{\gamma }{\alpha}\sum\limits_{i \in \mathcal{R}} {{{\left| {{\mathcal{B}_i}} \right|}^2}}+8L\alpha^2t_k + 4\left( {2L - \mu } \right){\alpha ^2}{D_F}\left( {{x_k},{x^*}} \right) \\
&- \frac{{c}}{{{q_{\max }}}}t_k + \frac{{c}}{{{q_{\min }}}}{D_F}\left( {{x_k},{x^*}} \right)\\
\le & \left( {1 - \left( {\gamma \alpha  - \frac{L}{2}\left( {4\left( {2L - \mu } \right){\alpha ^2} + \frac{c}{{{q_{\min }}}}} \right)} \right)} \right)\left\| {{x_k} - {x^*}} \right\|_2^2\\
&+ P_1{\alpha ^2} + {P_2}\alpha  + \left( {8L{\alpha ^2} - \frac{c}{{{q_{\max }}}}} \right)t_k,
\end{aligned}
\end{equation}
where the last inequality employs $L$-smoothness of the local objective function $f_i$, $\forall i \in \mathcal{R}$.
We proceed to choose $0 < \alpha  \le \gamma /\left( {8L\left( {2L - \mu } \right)} \right)$ and set ${{c} = \tilde c{\alpha }}$ with $0 < \tilde c \le {q_{\min }}\gamma /L$, such that (\ref{E10-14}) becomes
\begin{equation}\label{E10-15}
\begin{aligned}
&{\mathbb{E}_k}\left[ {\left\| {{x_{k + 1}} - {x^*}} \right\|_2^2} \right] + \frac{{{q_{\min }}\gamma \alpha }}{L}\left( {{\mathbb{E}_k}\left[ {t_{k + 1}} \right] - t_k} \right)\\
\le & \left( {1 \!-\! \frac{\gamma }{4}{\alpha }} \right)\left\| {{x_k} \!-\! {x^*}} \right\|_2^2 \!+\! \left( {8L{\alpha ^2} \!-\! \frac{c}{{{q_{\max }}}}} \right)t_k + P_1{\alpha ^2} \!+\! P_2\alpha.
\end{aligned}
\end{equation}
To proceed, via fixing $\tilde c = {q_{\min }}\gamma /L$ and $0 < \alpha  \le \gamma /\left( {8L\left( {2L - \mu } \right)} \right)$, it is equivalent to write (\ref{E10-15}) as
\begin{equation}\label{E10-16}
\begin{aligned}
&{\mathbb{E}_k}\left[ {\left\| {{x_{k + 1}} - {x^*}} \right\|_2^2} \right] + \frac{{{q_{\min }}\gamma {\alpha }}}{L}\left( {{\mathbb{E}_k}\left[ {t_{k + 1}} \right] - t_k} \right)\\
\le & \left( {1 - \frac{\gamma }{4}{\alpha}} \right)\left\| {{x_k} - {x^*}} \right\|_2^2 + 4n{{\phi }^2}\alpha ^2( {2\sum\limits_{i \in \mathcal{R}} {{{\left| {{\mathcal{R}_i}} \right|}^2}}  + \sum\limits_{i \in \mathcal{R}} {{{\left| {{\mathcal{B}_i}} \right|}^2}} } )\\
&+ \frac{{n{{\phi }^2}}}{\gamma }{\alpha }\sum\limits_{i \in \mathcal{R}} {{{\left| {{\mathcal{B}_i}} \right|}^2}} + \left( {8L{\alpha } - \frac{{\gamma{q_{\min }} }}{{L{q_{\max }}}}} \right){\alpha }t_k.
\end{aligned}
\end{equation}
We continue to define ${U_k}: = \left\| {{x_k} - {x^*}} \right\|_2^2 + {q_{\min }}\gamma {\alpha}t_k/ {L}$, which is non-negative due to $t_k \ge 0$. Based on this definition, if we select the constant step-size $0 < \alpha  \le 4\gamma /\left( {{\kappa _q}\left( {32{L^2} + {q_{\min }}{\gamma ^2}} \right)} \right)$, then it is natural to convert (\ref{E10-16}) into
\begin{equation}\label{E10-17}
\begin{aligned}
{\mathbb{E}_k}\left[ {{U_{k + 1}}} \right] \le & {\mathbb{E}_k}\left[ {\left\| {{x_{k + 1}} - {x^*}} \right\|_2^2} \right] + \frac{{{q_{\min }}\gamma \alpha }}{L}{\mathbb{E}_k}\left[ {t_{k + 1}} \right]\\
\le &\left( {1 - \frac{\gamma }{4}\alpha } \right)\left\| {{x_k} - {x^*}} \right\|_2^2 \!+\! \left( {1 - \frac{\gamma }{4}\alpha } \right)\frac{{{q_{\min }}\gamma \alpha }}{L}t_k\\
& + P_1{\alpha ^2} + {P_2}\alpha\\
= & \left( {1 - \frac{\gamma }{4}{\alpha }} \right){U_k} + P_1{\alpha ^2}+{P_2}\alpha.
\end{aligned}
\end{equation}
Summarizing all the upper bounds on the constant step-size generates a sufficient selection range as follows:
\begin{equation}\label{E10-18}
0 < \alpha  \le \frac{1}{{{\kappa _q}\left( {32{{\left( {1 + {\kappa _f}} \right)}^2} + {q_{\min }}} \right)}}\frac{1}{\mu }.
\end{equation}
Under the condition of (\ref{E10-18}), taking the total expectation on the both sides of (\ref{E10-17}) obtains
\begin{equation}\label{E10-19}
\mathbb{E}\left[ {{U_{k + 1}}} \right] \le \left( {1 - \frac{\gamma }{4}\alpha } \right)\mathbb{E}\left[ {{U_k}} \right] + {\alpha ^2}P_1 + \alpha P_2.
\end{equation}
Applying telescopic cancellation to (\ref{E10-19}) for $ k \ge 0$ obtains
\begin{equation}\label{E10-20}
\begin{aligned}
{\mathbb{E}}\left[ {\left\| {{x_{k + 1}} - {x^*}} \right\|_2^2} \right]\le & 4\left( {\frac{{{P_1}}}{\gamma }\alpha  + E} \right)\left( {1 - {{\left( {1 - \frac{\gamma }{4}\alpha } \right)}^{k+1}}} \right)\\
& + {\left( {1 - \frac{\gamma }{4}\alpha } \right)^{k+1}}{U_0},
\end{aligned}
\end{equation}
where ${U_0} = \left\| {{x_0} - {x^*}} \right\|_2^2 + {q_{\min }}\gamma \alpha {t_0}/L$. It is worthwhile to mention that by specifying $r_k$ and $t_k$ as $r^u_k$ and $t^u_k$ (resp., $r^w_k$ and $t^w_k$), the linear convergence rate is established for \textit{Prox-DBRO-SAGA} (resp., \textit{Prox-DBRO-LSVRG}).

\subsection{Proof of Theorem \ref{T3}}\label{Appen6}
Following the optimal condition given by (\ref{E10-1}), we have
\begin{equation}\label{E11-1}
{x^*} = {\mathbf{pro}}{{\mathbf{x}}_{{\alpha _k},G}}\left\{ {{x^*} - {\alpha _k}\left( {\nabla F\left( {{x^*}} \right) + {\partial _x}\chi \left( {{x^*}} \right)} \right)} \right\}
\end{equation}

We next proceed with the convergence analysis based on the transformed version (\ref{E11-1}) of the compact updates (\ref{E3-5-1})-(\ref{E3-5-2}).
\begin{equation*}
\begin{aligned}
&{\mathbb{E}_k}\left[ {\left\| {{x_{k + 1}} - {x^*}} \right\|_2^2} \right]\\
= & {\mathbb{E}_k}\left[ {\left\| {{\mathbf{pro}}{{\mathbf{x}}_{{\alpha _k},G}}\left\{ {{{\bar x}_k}} \right\} - {\mathbf{pro}}{{\mathbf{x}}_{{\alpha _k},G}}\left\{ {{x^*}} \right.} \right. - {\alpha _k}\nabla F\left( {{x^*}} \right)} \right.\\
& \left. {\left. { - {\alpha _k}{\partial _x}\chi \left( {{x^*}} \right)} \right\|_2^2} \right]\\
\le & {\mathbb{E}_k}\left[ {\left\| {{x_k} - {\alpha _k}\left( {r_k + {\partial _x}\chi \left( {{x_k}} \right) + {\partial _x}\delta \left( {{x_k}} \right)} \right)} \right.} \right.\\
&\left. { - \left. {\left( {{x^*} - {\alpha _k}\left( {\nabla F\left( {{x^*}} \right) + {\partial _x}\chi \left( {{x^*}} \right)} \right)} \right)} \right\|_2^2} \right]\\
= & \left\| {{x_k} - {x^*}} \right\|_2^2 - 2{\alpha _k}{\mathbb{E}_k}\left[ {\left\langle {{x_k} - {x^*},r_k - \nabla F\left( {{x^*}} \right)} \right\rangle } \right]\\
\end{aligned}
\end{equation*}
\begin{equation}\label{E11-2}
\begin{aligned}
& - 2{\alpha _k}\left\langle {{x_k} - {x^*},{\partial _x}\chi \left( {{x_k}} \right) - {\partial _x}\chi \left( {{x^*}} \right)} \right\rangle  \\
& - 2{\alpha _k}\left\langle {{x_k} - {x^*},{\partial _x}\delta \left( {{x_k}} \right)} \right\rangle + \alpha _k^2{\mathbb{E}_k}\left[ {\left\| {r_k + {\partial _x}\chi \left( {{x_k}} \right) } \right.} \right.\\
&\left. {\left. {+ {\partial _x}\delta \left( {{x_k}} \right) - \left( {\nabla F\left( {{x^*}} \right) + {\partial _x}\chi \left( {{x^*}} \right)} \right)} \right\|_2^2} \right],
\end{aligned}
\end{equation}
where the inequality applies the non-expansiveness of the proximal operator ${\mathbf{Pro}}{{\mathbf{x}}_{{\alpha _k},G}}\left\{ {\cdot} \right\}$.
To proceed, we continue to seek an upper bound for the last term in the RHS of (\ref{E11-2}) as follows:
\begin{equation}\label{E11-3}
\begin{aligned}
&{\mathbb{E}_k}\left[ {\left\| {r_k + {\partial _x}\chi \left( {{x_k}} \right) + {\partial _x}\delta \left( {{x_k}} \right) - \left( {\nabla F\left( {{x^*}} \right) + {\partial _x}\chi \left( {{x^*}} \right)} \right)} \right\|_2^2} \right]\\
\le & 4{\mathbb{E}_k}\left[ {\left\| {r_k - \nabla F\left( {{x^*}} \right)} \right\|_2^2} \right]  + 2\left\| {{\partial _x}\chi \left( {{x_k}} \right) - {\partial _x}\chi \left( {{x^*}} \right)} \right\|_2^2 \\
& + 4\left\| {{\partial _x}\delta \left( {{x_k}} \right)} \right\|_2^2\\
\le & 4\left( {4Lt_k + 2\left( {2L - \mu } \right){D_F}\left( {{x_k},{x^*}} \right)} \right) + 4\left\| {{\partial _x}\delta \left( {{x_k}} \right)} \right\|_2^2 \\
&+ 2\left\| {{\partial _x}\chi \left( {{x_k}} \right) - {\partial _x}\chi \left( {{x^*}} \right)} \right\|_2^2\\
\le & 4\left( {2Lt_k + \left( {2L - \mu } \right){D_F}\left( {{x_k},{x^*}} \right)} \right) + 8n{{\phi }^2}\sum\limits_{i \in \mathcal{R}} {{{\left| {{\mathcal{R}_i}} \right|}^2}} \\
& + 4n{{\phi }^2}\sum\limits_{i \in \mathcal{R}} {{{\left| {{\mathcal{B}_i}} \right|}^2}}.
\end{aligned}
\end{equation}
where the first inequality uses the basic inequality twice, the second inequality applies the result in Lemma \ref{L2}, and the last inequality is owing to (\ref{E10-4}) and (\ref{E10-6}). We proceed to handle the second term in the RHS of (\ref{E11-2}) as follows:
\begin{equation}\label{E11-4}
\begin{aligned}
&- 2{\mathbb{E}_k}\left[ {\left\langle {{x_k} - {x^*},r_k - \nabla F\left( {{x^*}} \right)} \right\rangle } \right]\\
= & - 2{\mathbb{E}_k}\left[ {\left\langle {{x_k} - {x^*},\nabla F\left( {{x_k}} \right) - \nabla F\left( {{x^*}} \right)} \right\rangle } \right] \\
\le &  \frac{{ - 2\mu L}}{{\mu  + L}}\left\| {{x_k} - {x^*}} \right\|_2^2 - \frac{2}{{\mu  + L}}\left\| {\nabla F\left( {{x_k}} \right) - \nabla F\left( {{x^*}} \right)} \right\|_2^2.
\end{aligned}
\end{equation}
Since $\chi $ is a convex function, the third term in the RHS of (\ref{E11-2}) can be dropped as the sequel inequality holds true
\begin{equation}\label{E11-5}
- \left\langle {{x_k} - {x^*},{\partial _x}\chi \left( {{x_k}} \right) - {\partial _x}\chi \left( {{x^*}} \right)} \right\rangle  \le 0. \\
\end{equation}
For any $\gamma > 0$, the fourth term in the RHS of (\ref{E11-2}) can be tackled similar to (\ref{E10-10}). Recalling the definitions of $ P_1$ and $P_2$ and plugging the results (\ref{E10-10}), (\ref{E11-3}), and (\ref{E11-5}) into (\ref{E11-2}) reduces to
\begin{equation}\label{E11-6}
\begin{aligned}
&{\mathbb{E}_k}\left[ {\left\| {{x_{k + 1}} - {x^*}} \right\|_2^2} \right]\\
\le & \left( {1 - \gamma{\alpha _k}} \right)\left\| {{x_k} - {x^*}} \right\|_2^2 + 8L\alpha _k^2t_k + P_1\alpha _k^2 + P_2{\alpha _k} \\
 &+ 4\left( {2L - \mu } \right)\alpha _k^2{D_F}\left( {{x_k},{x^*}} \right).
\end{aligned}
\end{equation}
According to Lemma \ref{L1}, we introduce an iteration-shifting variable ${c_k}>0$, such that
\begin{equation}\label{E11-7}
{c_k}\left( {{\mathbb{E}_k}\left[ {t_{k + 1}} \right] - t_k} \right) \le  - \frac{{{c_k}}}{{{q_{\max }}}}t_k + \frac{{{c_k}}}{{{q_{\min }}}}{D_F}\left( {{x_k},{x^*}} \right).
\end{equation}
Combining (\ref{E11-6}) and (\ref{E11-7}) obtains
\begin{equation*}
\begin{aligned}
&{\mathbb{E}_k}\left[ {\left\| {{x_{k + 1}} - {x^*}} \right\|_2^2} \right] + {c_k}\left( {{\mathbb{E}_k}\left[ {t_{k + 1}} \right] - t_k} \right)\\
\le & \left( {1 - \gamma {\alpha _k}} \right)\left\| {{x_k} - {x^*}} \right\|_2^2 + \left( {8L\alpha _k^2 - \frac{{{c_k}}}{{{q_{\max }}}}} \right)t_k + P_1\alpha _k^2\\
&+P_2{\alpha _k} + \left( {\frac{{{c_k}}}{{{q_{\min }}}}+4\left( {2L - \mu } \right)\alpha _k^2} \right){D_F}\left( {{x_k},{x^*}} \right)\\
\le & \left( {1 - \left( {\gamma {\alpha _k} - \frac{L}{2}\left( {4\left( {2L - \mu } \right)\alpha _k^2 + \frac{{{c_k}}}{{{q_{\min }}}}} \right)} \right)} \right)\left\| {{x_k} - {x^*}} \right\|_2^2\\
\end{aligned}
\end{equation*}
\begin{equation}\label{E11-8}
\begin{aligned}
&+ \left( {8L\alpha _k^2 - \frac{{{c_k}}}{{{q_{\max }}}}} \right)t_k +P_1\alpha _k^2 + {P_2}{\alpha _k},
\end{aligned}
\end{equation}
where the last inequality uses $L$-smoothness of the local objective function $f_i$, $\forall i \in \mathcal{R}$. Via setting ${{c_k} = \tilde c{\alpha _k}}$ and $0 < {\alpha _k} \le \gamma /\left( {8L\left( {2L - \mu } \right)} \right)$, we have
\begin{equation}\label{E11-9}
\begin{aligned}
& {\mathbb{E}_k}\left[ {\left\| {{x_{k + 1}} - {x^*}} \right\|_2^2} \right] + \frac{{\gamma {q_{\min }}{\alpha _k}}}{L}{\mathbb{E}_k}\left[ {t_{k + 1}} \right]\\
\le & \left( {1 - \frac{\gamma }{4}} \right)\left\| {{x_k} - {x^*}} \right\|_2^2 + \left( {( {1 - \frac{1}{{{q_{\max }}}}} )\frac{{{q_{\min }}\gamma }}{L} + 8L{\alpha _k}} \right){\alpha _k}t_k\\
& + P_1\alpha _k^2 + {P_2}{\alpha _k}.
\end{aligned}
\end{equation}
We define ${\tilde U_k}: = \left\| {{x_k} - {x^*}} \right\|_2^2 + {q_{\min }}\gamma {\alpha _k}t_k/L$, which is non-negative, since $t_k$ is non-negative. We further set $0 < {\alpha _k} \le 4\gamma /\left( {{\kappa _q}\left( {32{L^2} + {q_{\min }}{\gamma ^2}} \right)} \right)$ and take the total expectation on the both sides of (\ref{E11-9}) to obtain
\begin{equation}\label{E11-10}
\begin{aligned}
{\mathbb{E}}\left[ {{\tilde U_{k + 1}}} \right]\le & {\mathbb{E}}\left[ {\left\| {{x_{k + 1}} - {x^*}} \right\|_2^2} \right] + \frac{{{q_{\min }}\gamma {\alpha _k}}}{L}{\mathbb{E}}\left[ {t_{k + 1}} \right]\\
\le & \left( {1 - \frac{\gamma }{4}{\alpha _k}} \right){\mathbb{E}}\left[{\tilde U_k}\right] + P_1\alpha _k^2 + {P_2}{\alpha _k},
\end{aligned}
\end{equation}
where the first inequality is due to the fact that the step-size $\alpha_k$ is decaying. By summarizing all the required upper bounds on the step-size, it suffices to consider a feasible range as follows:
\begin{equation}\label{E11-11}
0 < \alpha_k  \le \frac{1}{{{\kappa _q}\left( {32{{\left( {1 + {\kappa _f}} \right)}^2} + {q_{\min }}} \right)}}\frac{1}{\mu }.
\end{equation}
According to (\ref{E11-11}), we set ${\alpha _k} = \theta /\left( {k + \xi} \right)$, $\forall k \ge 0$, with $\theta  > 4/\gamma $ and $\xi  = {\kappa _q}\left( {32{{\left( {1 + {\kappa _f}} \right)}^2} + {q_{\min }}} \right)\mu \theta $. We next prove
\begin{equation}\label{E11-12}
\mathbb{E}\left[ {{\tilde U_k}} \right] \le \Xi /\left( {k + \xi} \right) + \tilde E, \forall k \ge 0,
\end{equation}
by induction. Firstly, for $k = 0$, we know that
\begin{equation}\label{E11-13}
{\tilde U_1} \le \left( {1 - \frac{\gamma }{4}{\alpha _0}} \right){\tilde U_0} + \alpha _0^2{P_1} + {\alpha _0}{P_2}.
\end{equation}
Therefore, for a bounded and positive constant $\tilde E$, if $\Xi  \ge \left( {\xi  - \gamma\theta/4 } \right){\tilde U_0} + {{\theta ^2}{P_1}/\xi }  + \theta {P_2} - \xi \tilde E$, we have
\begin{equation}\label{E11-14}
\tilde U_1 \le \left( {1 - \frac{\gamma }{4}{\alpha _0}} \right){\tilde U_0} + \alpha _0^2{P_1} + {\alpha _0}{P_2} \le \frac{\Xi }{\xi } + \tilde E,
\end{equation}
with ${\alpha _0} = \theta /\xi $. We assume that for $k = K$, $ K \ge 1$, it satisfies that
\begin{equation}\label{E11-15}
\begin{aligned}
\mathbb{E}\left[ {{\tilde U_{K + 1}}} \right] \le &\left( {1 - \frac{\gamma }{4}{\alpha _{K}}} \right)\mathbb{E}\left[ {{\tilde U_{K}}} \right] + \alpha _{K}^2{P_1} + {\alpha _{K}}{P_2}\\
\le & \frac{\Xi }{{K + \xi}} + \tilde E.
\end{aligned}
\end{equation}
Then, we will prove that for $k = K + 1$,
\begin{equation}\label{E11-16}
\mathbb{E}\left[ {{\tilde U_{K + 2}}} \right] \le  \frac{\Xi }{{K + \xi + 1}} + \tilde E,
\end{equation}
holds true. We define $\tilde \gamma  := \gamma /4$ and set $\tilde E \ge {P_2}/\tilde \gamma $ and $\Xi  \ge {\theta ^2}{P_1}/\left( {\tilde \gamma \theta  - 1} \right)$ with $\theta  > 1/\tilde \gamma $. We have
\begin{equation*}
\begin{aligned}
& \mathbb{E}\left[ {{\tilde U_{K + 2}}} \right] \\
\le & \left( {1 - \tilde \gamma {\alpha _{K + 1}}} \right)\mathbb{E}\left[ {{\tilde U_{K + 1}}} \right] + \alpha _{K + 1}^2{P_1} + {\alpha _{K + 1}}{P_2}\\
\le & \left( {1 - \frac{{\tilde \gamma \theta }}{{K + \xi + 1}}} \right)\left( {\frac{\Xi }{{K + \xi}} + \tilde E} \right) + \frac{{{\theta ^2}}}{{{{\left( {K + \xi + 1} \right)}^2}}}{P_1}\\
& + \frac{\theta }{{K + \xi + 1}}{P_2}\\
\end{aligned}
\end{equation*}
\begin{equation}\label{E11-17}
\begin{aligned}
\le &\left( {1 - \frac{{\tilde \gamma \theta }}{{K + \xi + 1}}} \right)\frac{\Xi }{{K + \xi}} + \tilde E + \frac{{{\theta ^2}}}{{{{\left( {K + \xi + 1} \right)}^2}}}{P_1}\\
\le & \left( {1 - \frac{{\tilde \gamma \theta }}{{K + \xi + 1}}} \right)\frac{\Xi }{{K + \xi}} + \tilde E + \frac{{\Xi \left( {\tilde \gamma \theta  - 1} \right)}}{{{{\left( {K + \xi + 1} \right)}^2}}}\\
\le & \left( {1 - \frac{{\tilde \gamma \theta }}{{K + \xi + 1}}} \right)\frac{\Xi }{{K + \xi}} + \frac{{\Xi \left( {\tilde \gamma \theta  - 1} \right)}}{{\left( {K + \xi + 1} \right)\left( {K + \xi} \right)}} + \tilde E \\
= &  \frac{\Xi }{{K + \xi + 1}} + \tilde E,
\end{aligned}
\end{equation}
which means the relation (\ref{E11-12}) holds true. Via replacing $\tilde E$ with its upper bound $E = {P_2}/\tilde \gamma $, it is straightforward to verify
\begin{equation}\label{E11-18}
{\mathbb{E}}\left[ {\left\| {{x_{k}} - {x^*}} \right\|_2^2} \right] \le  \frac{\Xi }{{k  + \xi }} + E, \forall k \ge 0,
\end{equation}
owing to $t_k \ge 0$. Through specifying $r_k$ and $t_k$ as $r^u_k$ and $t^u_k$ (resp., $r^w_k$ and $t^w_k$), the sub-linear convergence rate is established for \textit{Prox-DBRO-SAGA} (resp., \textit{Prox-DBRO-LSVRG}).

\bibliographystyle{IEEEtran}
\bibliography{Prox-DBRO-VR}
\end{document}